%% file: pt2.tex
\documentclass[12pt]{amsart}
\usepackage[dvips]{graphicx,color}
\usepackage[dvips]{epsfig}
\usepackage{amsmath,amssymb,latexsym,amscd,xypic}
\DeclareRobustCommand{\SkipTocEntry}[4]{}
\makeatletter
\newcommand\@dotsep{4.5}
\def\@tocline#1#2#3#4#5#6#7{\relax
  \ifnum #1>\c@tocdepth % then omit
  \else
    \par \addpenalty\@secpenalty\addvspace{#2}%
    \begingroup \hyphenpenalty\@M
    \@ifempty{#4}{%
      \@tempdima\csname r@tocindent\number#1\endcsname\relax
    }{%
      \@tempdima#4\relax
    }%
    \parindent\z@ \leftskip#3\relax \advance\leftskip\@tempdima\relax
    \rightskip\@pnumwidth plus1em \parfillskip-\@pnumwidth
    #5\leavevmode\hskip-\@tempdima #6\relax
    \leaders\hbox{$\m@th
      \mkern \@dotsep mu\hbox{.}\mkern \@dotsep mu$}\hfill
    \hbox to\@pnumwidth{\@tocpagenum{#7}}\par
    \nobreak
    \endgroup
  \fi}
\makeatother 

\DeclareFontFamily{OT1}{rsfs}{}
\DeclareFontShape{OT1}{rsfs}{n}{it}{<-> rsfs10}{}
\DeclareMathAlphabet{\curly}{OT1}{rsfs}{n}{it}

\newcommand\mmuu{{\stackrel{\to}{\mu}}}
\newcommand{\cO}{\mathcal{O}}
\newcommand{\N}{\mathcal{N}}

\newcommand{\T}{\mathbf{T}}
\newcommand\C{\mathbb C}
\newcommand\F{\curly F}
\newcommand{\Ib}{{\mathbb{I}}}
\newcommand\I{\curly I}
\newcommand\II{\mathrm{I\hspace{-1.3pt}I}}
\newcommand\III{\mathrm{I\hspace{-1.3pt}I\hspace{-1.3pt}I}}

\newcommand\m{\mathfrak m}

\renewcommand\O{\mathcal O}
\newcommand\PP{\mathbb P}
\newcommand\FF{\mathbb F}
\newcommand\FFF{\mathsf{F}}
\newcommand\GGG{\mathsf{G}}

\newcommand\Q{\mathbb Q}

\newcommand\Z{\mathbb Z}
\newcommand\bW{\mathsf{W}}

\newcommand{\bV}{\mathsf{V}}
\newcommand{\bE}{\mathsf{E}}
\newcommand{\bw}{\mathsf{w}}

\newcommand{\Rt}[1]{\stackrel{#1\,}{\longrightarrow}}
\newcommand\To{\longrightarrow}
\newcommand\into{\hookrightarrow}
\newcommand\Into{\ar@{^{ (}->}[r]}

\renewcommand\_{^{}_}

\newcommand\bull{{\scriptscriptstyle\bullet}}
\newcommand\udot{^\bull}

\newcommand\tr{\operatorname{tr}}

\newcommand\Hom{\operatorname{Hom}}
\renewcommand\hom{\curly H\!om}
\newcommand\Ext{\operatorname{Ext}}

\newcommand{\extt}{\operatorname{Ext}}
\newcommand{\sheafext}{\mathcal{E}xt}

\newcommand\beq[1]{\begin{equation}\label{#1}}
\newcommand\eeq{\end{equation}}
\newcommand\beqa{\begin{eqnarray*}}
\newcommand\eeqa{\end{eqnarray*}}

\makeatletter \@addtoreset{equation}{section} \makeatother

\newtheorem{thm}{Theorem}
\newtheorem{thmc}{Theorem/Conjecture}
\newtheorem{lem}{Lemma}

\newtheorem{conj}{Conjecture}
\newtheorem{prop}{Proposition}

\title{\textbf{The $3$-fold vertex via stable pairs}}
\author{R. Pandharipande and R. P. Thomas}
\date{September 2007}

\begin{document}

\begin{abstract} \noindent
The theory of stable pairs in the derived category yields
an enumerative geometry of curves in $3$-folds.
We evaluate the equivariant  
vertex for stable pairs on toric $3$-folds
in terms of weighted box counting.
In the toric Calabi-Yau case, the result simplifies to a new
form of pure box counting. 
The conjectural equivalence
with the DT vertex predicts remarkable identities.

The equivariant vertex governs primary insertions
in the theory of stable pairs for toric varieties.
We consider also the descendent vertex and
conjecture the complete rationality of the
descendent theory for stable pairs.

\end{abstract}

\maketitle

\setcounter{tocdepth}{1} 
\tableofcontents

%%%%%%%%%%%%%%%%%%%%%%%%%%%%%%%%%%%%%%%%%%%%%%%%%%%%%%%%%%%%%%%%%%%%%%%%%%%

\setcounter{section}{-1}
\section{Introduction}

\subsection{Overview}
Let $X$ be a nonsingular 3-fold, and let
$$\beta \in H_2(X,\mathbb{Z})$$ be a non-zero class. We are interested here
in the moduli space of stable pairs
$$[\O_X \stackrel{s}{\rightarrow} F] \in P_n(X,\beta)$$
where $F$ is a pure sheaf supported on a Cohen-Macaulay subcurve of $X$, 
$s$ is a morphism with 0-dimensional cokernel, and
$$\chi(F)=n, \  \  \ [F]=\beta.$$
The space $P_n(X,\beta)$
carries a virtual fundamental class obtained from the 
deformation theory of complexes in
the derived category \cite{pt}. A review can be found in Section \ref{ooo}.

If $X$ is toric, we may calculate the stable pairs invariants by localization
with respect to the torus action \cite{GraberP}. 
The outcome is expressed in terms of
the associated polyhedron $\Delta(X)$. 
The edge
contributions
are related to partition sums. 
The vertex contributions, related to box counting,
are the most interesting aspect of the geometry.

We calculate the edge and vertex contributions for toric $X$ in terms
of weighted partition and box counts. In case $X$ is toric Calabi-Yau, the
formulas simplify to pure box counting. 
The subject is related to dualities in string theory, wall-crossing
formulae in the derived category, commutative algebra, and the
combinatorics of 3-dimensional partitions.

\subsection{Toric geometry}
 \label{tgeom}
Let $X$ be a nonsingular toric 3-fold acted upon by a
3-dimensional complex torus $\T$.
Let $\Delta(X)$ denote the
Newton polyhedron of $X$ determined by a polarization.
The polyhedron $\Delta(X)$ is the image of $X$ under
the moment map.

The vertices of the polyhedron $\Delta(X)$ correspond to fixed points
$$X^\T=\{X_\alpha\}$$
of the $\T$-action. 
For each $X_\alpha$, there is a canonical,
$\T$-invariant, affine open chart,
$$U_\alpha\cong \C^3,$$ centered
at $X_\alpha$. We may choose coordinates $t_i$ on
$\T$ and coordinates $x_i$ on
$U_\alpha$ for which
the $\T$-action on $U_\alpha$
is
determined by 
\begin{equation}
  \label{st_act}
  (t_1,t_2,t_3)\cdot x_i = t_i x_i  \,.
\end{equation}
%%
%In these coordinates,
%the tangent representation $X_\alpha$
% has character $$t^{-1}_1+ t_2^{-1}+ t_3^{-1}.$$
%We will use the covering $\{ U_\alpha\}$ of $X$ to compute
%Cech cohomology.

The edges of $\Delta(X)$ correspond to the 
 $\T$-invariant lines of $X$.
More precisely, if
$$C_{\alpha\beta}\subset  X$$
is a $\T$-invariant line incident to the fixed points
$X_\alpha$ and $X_\beta$, then
$C_{\alpha\beta}$ corresponds to an edge of $\Delta(X)$
joining the vertices $X_\alpha$ and $X_\beta$.

The geometry of $\Delta(X)$ near the edge is determined by the normal
bundle $\N_{C_{\alpha\beta}/X}$. If
$$
\N_{C_{\alpha\beta}/X} = \cO(m_{\alpha\beta}) \oplus \cO(m'_{\alpha\beta})
$$
then the transition functions between the charts $U_\alpha$ and
$U_\beta$ can be taken to be of the form
\begin{equation}\label{trans}
  (x_1,x_2,x_3) \mapsto (x_1^{-1}, x_2 \,x_1^{-m_{\alpha\beta}} ,
x_3 \,x_1^{-m'_{\alpha\beta}} ) \,.
\end{equation}
The curve $C_{\alpha\beta}$ is then 
defined in these coordinates by $x_2=x_3=0$.
If $X$ is Calabi-Yau, then degree of $\N_{C_{\alpha\beta}/X}$ is $-2$ and
$$m_{\alpha\beta} = m'_{\alpha\beta} \mod 2.$$

\subsection{Localization}
The $\T$-action on $X$  
canonically induces a $\T$-action on the moduli space of
pairs $P_n(X,\beta)$. Our first result is a determination of the
$\T$-fixed loci of $P_n(X,\beta)$. Let
$$\mathbf{Q} \subset P_n(X,\beta)^\T$$
be a connected $\T$-fixed locus.

\begin{thm} \label{ccc} $\mathbf{Q}$ is a product of $\PP^1$s.
\end{thm}

The $0^{th}$ product of $\PP^1$ is a point. Indeed, if $X$ is a local 
toric
surface, the $\T$-fixed points of $P_n(X,\beta)$ are isolated. 
Positive dimensional fixed loci occur only in the fully 3-dimensional
setting. 
Theorem \ref{ccc} is proven by an explicit characterization of 
the $\T$-fixed points in terms of box configurations in Sections 
\ref{ttt} and \ref{tang}.

Let $[\O_X \rightarrow F] \in \mathbf{Q}$
be a stable pair, and let
$C\subset X$
be the Cohen-Macaulay curve obtained from the (scheme-theoretic) 
support of $F$. Certainly $C$ has set-theoretic
support on the edge curves
$$\bigcup_{\alpha,\beta} C_{\alpha\beta}\subset X.$$
Since $C$ must be $\T$-invariant, $C$ determines 
a partition $\mu_{\alpha\beta}$
at each edge. The size $|\mu_{\alpha\beta}|$ of the partition 
is simply the multiplicity of $C$ along $C_{\alpha\beta}$.
The partition $\mu_{\alpha\beta}$ is the same for each stable
pair in $\mathbf{Q}$. All the moduli in $\mathbf{Q}$ are obtained from the
vertices.

A complete determination of the $\T$-equivariant
contribution of the $\T$-fixed locus $\mathbf{Q}$ to the
stable pairs theory of $X$ is the main calculation of our paper.
The result is easiest to state in the toric Calabi-Yau case for the
basic stable pairs invariant 
\begin{equation}\label{bdww}
P_{n,\beta} = \int_{[P_n(X,\beta)]^{vir}} 1.
\end{equation}
%In the toric Calabi-Yau case, let
%$$\T_{CY} \subset \T$$
%be the 2-dimensional torus which preserves a non-vanishing
%canonical section.
%The integral \eqref{bdww} is defined equivariantly with
%respect to $\T_{CY}$.
If $X$ is toric Calabi-Yau, define 
the restricted contribution of 
$\mathbf{Q}$
to the invariant 
$P_{n,\beta}$
by 
$$\chi\_{top}(\mathbf{Q})\cdot  (-1)^{n + \sum_{\alpha \beta}m_{\alpha\beta}
|\mu_{\alpha\beta}|}$$
where $\chi\_{top}$ is the topological Euler characteristic.

\setcounter{thmc}{1}
\begin{thmc} \label{ddd} 
The toric Calabi-Yau
 invariant $P_{n,\beta}$ is obtained by summing over all components
of $P_n(X,\beta)^\T$, 
\begin{eqnarray*}
Z_{P,\beta}(q) & =&  \sum_{n} P_{n,\beta}\ q^n\\
&  =& \sum_{n} \sum_{\mathbf{Q} \subset P_n(X,\beta)^\T}          
\chi\_{top}(\mathbf{Q})\cdot  (-1)^{n + \sum_{\alpha \beta} m_{\alpha\beta}
|\mu^{\mathbf{Q}}_{\alpha\beta}|}\ q^n.
\end{eqnarray*}
\end{thmc}

We prove Theorem \ref{ddd} in the local Calabi-Yau
toric surface case (where
all vertices have at most 2 legs). For the 3-leg case, 
 our derivation
 at present 
depends upon conjectural{\footnote{The Calabi-Yau case is a particular
limit of the full $\T$-equivariant calculation. The conjectural
properties are needed to take the limit.}} properties of the
stable pair space, see Sections \ref{ccooj1}-\ref{ccooj2}. 
We will show the summation of Theorem \ref{ddd}
is a form of box counting.

\subsection{Correspondence with DT theory}
DT theory \cite{DT, Thomas} 
is defined by integration against the virtual fundamental
class of the moduli space $I_n(X,\beta)$ of ideal sheaves{\footnote{$I_n(X,\beta)$
is isomorphic to the Hilbert scheme.}}
$$0 \rightarrow \mathcal{I} \rightarrow \O_X \rightarrow \O_Y \rightarrow 0$$
satisfying
$$\chi(\O_Y)=n, \ \ \ [\O_Y]=\beta\in H_2(X,\mathbb{Z}).$$
In the Calabi-Yau case, the basic invariants are 
$$I_{n,\beta}= \int_{[I_n(X,\beta)]^{vir}} 1.$$

For  toric Calabi-Yau 3-folds,
the DT invariants have been calculated by
localization in \cite{MNOP1,MNOP2},
\begin{eqnarray}
\label{fred}
Z_{DT,\beta}(q) & =&  \sum_{n} I_{n,\beta}\ q^n\\ \nonumber
&  =& \sum_{n} \sum_{[\mathcal{I}] \subset I_n(X,\beta)^\T}          
(-1)^{n + \sum_{\alpha \beta} m_{\alpha\beta}
|\mu^{\mathcal{I}}_{\alpha\beta}|}q^n.
\end{eqnarray}

The result \eqref{fred} is parallel to Theorem \ref{ddd}.
The edge contributions in DT theory agree exactly with the
edge contributions in the theory of stable pairs. 
The main difference
occurs in the vertex contributions.
Since the fixed point set 
$$I_n(X,\beta)^\T\subset I_n(X,\beta)$$ consists
of isolated points, 
the DT result \eqref{fred} is easier to
prove than Theorem \ref{ttt} for stable pairs.
However, the stable pairs result is free of the
irrationalities related to unrestricted box counting.

The stable pairs theory is conjectured to be equivalent to DT theory
for all 3-folds \cite{pt}. In case $X$ is toric, the conjecture
specializes to the claim
$$Z_{P,{\beta}}(q) = \frac{Z_{{DT},\beta}(q)}{M(-q)^{|X^\T|}}\,,$$
where $M(-q)$ is the MacMahon function
$$M(-q) = \prod_{n\geq 1} \frac{1}{(1-(-q)^n)^n}\,.$$
The toric equivalence can be further refined to relate only
the vertex contributions of the two theories.{\footnote{
The refinement in the Calabi-Yau case is
Conjecture 5.1 of \cite{pt}.}}\ This is discussed in Section \ref{eqver}.

\subsection{Descendents}\label{des}
Let $X$ be a nonsingular projective 3-fold.
Let
$$\FF \rightarrow X\times P_{n}(X,\beta)$$
 denote the universal sheaf.{\footnote{The existence of the
universal sheaf is shown in Section 2.3 of \cite{pt}.}}
For a stable pair $[\O_X\to F]\in P_{n}(X,\beta)$, the restriction of
$\FF$
to the fiber
 $$X \times [\O_X \to F] \subset 
X\times P_{n}(X,\beta)
$$
is canonically isomorphic to $F$.
Let
$$\pi_X\colon X\times P_{n}(X,\beta)\to X,$$
$$\pi_P\colon X\times P_{n}(X,\beta)
\to P_{n}(X,\beta)$$
 be the projections onto the first and second factors.
Since $X$ is nonsingular
and
$\FF$ is $\pi_P$-flat, $\FF$ has a finite resolution 
by locally free sheaves.
Hence, the Chern character of the universal sheaf $\FF$ on 
$X \times P_n(X,\beta)$ is well-defined.
By definition, the operation
$$
\pi_{P*}\big(\pi_X^*(\gamma)\cdot \text{ch}_{2+i}(\FF)
\cap(\pi_P^*(\ \cdot\ )\big)\colon 
H_*(P_{n}(X,\beta))\to H_*(P_{n}(X,\beta))
$$
is the action of the descendent $\tau_i(\gamma)$, where
$\gamma \in H^*(X,\Z)$.

For nonzero $\beta\in H_2(X,\Z)$ and arbitrary $\gamma_i\in H^*(X,\Z)$,
define the stable pairs invariant with descendent insertions by
\begin{eqnarray*}
\left\langle \prod_{j=1}^k \tau_{i_j}(\gamma_j)
\right\rangle_{\!n,\beta}^{\!X}&  = &
\int_{[P_{n}(X,\beta)]^{vir}}
\prod_{j=1}^k \tau_{i_j}(\gamma_j) \\
& = & 
\int_{P_n(X,\beta)} \prod_{j=1}^k \tau_{i_j}(\gamma_{j})
\Big( [P_{n}(X,\beta)]^{vir}\Big).
\end{eqnarray*}
The partition function is 
$$
Z_{P,\beta}\left(   \prod_{j=1}^k \tau_{i_j}(\gamma_{k})
\right)
=\sum_{n} 
\left\langle \prod_{j=1}^k \tau_{i_j}(\gamma_{j}) 
\right\rangle_{\!n,\beta}^{\!X}q^n.
$$
Since $P_n(X,\beta)$ is empty for sufficiently negative
$n$, 
$Z_{P,\beta}\big(   \prod_{j=1}^k \tau_{i_j}(\gamma_{j})
\big)$
is a Laurent series in $q$.

\begin{conj}
\label{111} 
The partition function
$Z_{P,\beta}\big(   \prod_{j=1}^k \tau_{i_j}(\gamma_{j})
\big)$ is the 
Laurent expansion of a rational function in $q$.
\end{conj}

The partition functions with primary insertions (all $i_j=0$)
were conjectured to be rational and, furthermore, conjectured
to 
take a very restrictive BPS form in \cite{pt}. 
The analogue of 
BPS invariants in the presence of descendents is
an interesting question.

The descendent series of both Gromov-Witten theory and
DT theory are known to contain irrationalities.
Conjecture \ref{111} predicts the descendent theory of
stable pairs is much better behaved.

\subsection{Vertices}
The stable pairs vertices for toric 3-folds in increasing degree of
generality are:
\begin{enumerate}
\item[(i)] the toric Calabi-Yau vertex,
\item[(ii)] the equivariant vertex,
\item[(iii)] the equivariant descendent vertex.
\end{enumerate}
The vertices (i) and (ii) are discussed in Section \ref{eqver} and
\ref{cyver}.
We treat the localization formulas
for the descendent theory 
in Section \ref{desc}. 

\subsection*{Acknowledgements}  
We thank 
J. Bryan, E. Diaconescu, C. Faber,
D. Joyce, 
S. Katz, A. Klemm, D. Maclagan,
D. Maulik, G. Moore, A. Oblomkov, A. Okounkov, and
 S. Payne for conversations
related to stable pairs and box counting.
P. Hort helped with the  figures.

R.P. was partially supported by NSF grant DMS-0500187 and a Packard foundation
fellowship. 
 R.T. was partially supported
by a Royal Society University Research Fellowship.

R.T. would like to thank the Leverhulme Trust and Columbia University for a
visit to New York in the spring of 2007 when the project was started. 
Lectures on the results were given by R.P.
at the Centre de Recherches Math\'ematiques in Montr\'eal 
in the
summer of 2007.

\section{Stable pairs on $3$-folds}
\label{ooo}
\subsection{Definitions}
Let $X$ be a nonsingular quasi-projective $3$-fold over $\mathbb{C}$
with 
polarization $L$.
Let $\beta\in H_2(X,\mathbb{Z})$ be a nonzero class.
The moduli space $P_n(X,\beta)$ parameterizes \emph{stable pairs}
\begin{equation}\label{vqq2}
\O_X \stackrel{s}{\rightarrow} F
\end{equation}
where $F$ is a sheaf with Hilbert polynomial
$$ \chi(F\otimes L^k) = k\int_\beta c_1(L) + n$$
 and $s\in H^0(X,F)$ is a section. 
The two stability conditions are: 
\begin{enumerate}
\item[(i)]
the sheaf $F$ is {pure} with proper support, 
\item[(ii)] the section $\O_X \stackrel{s}{\rightarrow} F$ has 0-dimensional
cokernel.
\end{enumerate}
By definition, {\em purity} (i) means
every nonzero
subsheaf of $F$ has support of dimension 1 \cite{HLShaves}. In particular,
 purity implies  the (scheme theoretic) 
support $C_F$ of $F$ is a Cohen-Macaulay curve. 
A quasi-projective moduli space of stable pairs
 can be constructed by  a standard GIT analysis of Quot scheme
quotients \cite{LPPairs1}.

For convenience, we will often refer to the stable pair 
\eqref{vqq2} on $X$ simply by $(F,s)$.

\subsection{Virtual class}
A central result of \cite{pt} is the construction of a
virtual class on $P_n(X,\beta)$.
The standard approach to the deformation theory of pairs
fails to yield an appropriate 2-term deformation theory
for $P_n(X,\beta)$.
Instead, $P_n(X,\beta)$ is viewed in \cite{pt}
as a moduli space
of complexes in the derived category.

Let $D^b(X)$ be the bounded derived category of coherent
sheaves on $X$.
Let
$${I}\udot = \left\{ \O_X \rightarrow F \right\}\in
D^b(X)$$
be the complex determined by a stable pair.
The tangent-obstruction theory obtained by deforming ${I}\udot$
in $D^b(X)$ while fixing its determinant is 2-term and governed by the 
groups{\footnote{The subscript 0 denotes traceless $\Ext$.}}
$$\Ext^1({I}\udot, {I}\udot)_0, \ \
\Ext^2({I}\udot, {I}\udot)_0.$$
The virtual class 
$$[P_n(X,\beta)]^{vir} \in A_{\text{dim}^{vir}}
%{\int_\beta c_1(T_X)}
\left(P_n(X,\beta),\mathbb{Z}\right)$$
is then obtain by standard methods \cite{bf,lt}.
The virtual dimension is
$$\text{dim}^{vir} = \int_\beta c_1(T_X).$$

Apart from the derived category deformation theory,
the construction of the virtual class of $P_n(X,\beta)$
is parallel to virtual class construction in DT theory \cite{Thomas}.

\subsection{Characterization}
Consider the kernel/cokernel exact sequence associated to a stable
pair $(F,s)$,
\beq{IOFQ}
0\to\I_{C_F}\to\O_X\Rt{s}F\to Q\to0.
\eeq
The kernel is the ideal sheaf of the Cohen-Macaulay support
curve $C_F$ by Lemma 1.6 of \cite{pt}. The cokernel
$Q$ has dimension 0 support by stability.
The {\em reduced} support scheme, $\text{Support}^{red}(Q)$, is 
called the {\em zero locus} of the pair.
The zero locus lies on $C_F$.

Let $C\subset X$ be a fixed Cohen-Macaulay curve.
Stable pairs with support  $C$ and bounded zero locus are characterized
as follows.
Let $$\m\subset\O_C$$ 
be the ideal in $\O_C$ of a 0-dimensional subscheme. 
Since $$\hom(\m^r/\m^{r+1},\O_C)=0$$ by
the purity of $\O_C$, we obtain an inclusion $$\hom(\m^r,\O_C)\subset
\hom(\m^{r+1},\O_C).$$ 
The inclusion $\m^r\into\O_C$ induces a canonical section 
$$\O_C\into\hom(\m^r,\O_C).$$

\begin{prop} \label{descl}
A stable pair $(F,s)$ with support $C$ satisfying
$$\text{\em Support}^{red}(Q) \subset \text{\em Support}(\O_C/\m)$$ 
is equivalent to a subsheaf
of $\hom(\m^r,\O_C)/\O_C,\ r\gg0.$
\end{prop}

Alternatively, we may work with coherent subsheaves of the quasi-coherent sheaf
\begin{equation}\label{infhom}
\lim\limits_{\To}\hom(\m^r,\O_C)/\O_C
\end{equation}
Under the equivalence of Proposition \ref{descl}, the 
subsheaf of \eqref{infhom} corresponds to $Q$, giving a subsheaf $F$ of
$\lim\limits_{\To}\hom(\m^r,\O_C)$ containing the canonical subsheaf $\O_C$
and the sequence
$$0\to\O_C\stackrel{s}{\rightarrow} F \rightarrow Q \rightarrow 0.$$
 Proposition \ref{descl} is proven in \cite{pt}.

\section{$\T$-fixed points} \label{ttt}

\subsection{Affine charts}
Let $X$ be a nonsingular, quasi-projective, toric 3-fold, and let
\begin{equation}\label{vqaa}
[\O_X \stackrel{s}{\rightarrow} F] \in P_n(X,\beta)^\T
\end{equation}
be a $\T$-fixed stable pair.

Let $X_\alpha\in X^\T$ be a $\T$-fixed point with
associated $\T$-invariant affine chart
$U_\alpha\subset X$.
The restriction of the stable pair \eqref{vqaa} to $U_\alpha$,
\begin{equation}\label{vvvt}
\O_{U_\alpha} \stackrel{s_\alpha}{\rightarrow} F_\alpha
\end{equation}
determines an invariant section $s_\alpha$ of an
equivariant sheaf $F_\alpha$.

Let $x_1,x_2,x_3$ be coordinates on the affine chart $U_\alpha$
in which the $\T$-action takes the diagonal form.
$$(t_1,t_2,t_3) \cdot x_i = t_i x_i.$$
We will characterize the restricted data $(F_\alpha,s_\alpha)$
in the coordinates $x_i$.

\subsection{Monomial ideals and partitions}
Let $x_1,x_2$ be coordinates on the plane 
$\C^2$.
A subscheme $S\subset \C^2$ invariant under the 
action of the diagonal torus,
$$(t_1,t_2)\cdot x_i = t_ix_i$$
must be defined by a monomial ideal
$\I_S \subset \C[x_1,x_2]$.
If
$$\dim_\C \C[x_1,x_2]/\I_S < \infty$$
then $\I_S$ determines a finite partition $\mu_S$
by considering lattice points corresponding
to monomials of $\C[x_1,x_2]$
{\em not} contained in $\I_S$. 

Conversely, each partition $\mu$ determines a monomial ideal
$$\mu[x_1,x_2]\subset \C[x_1,x_2].$$
The ideal associated to the finite partition $(4,4,3,1,1)$ is 
displayed in Figure 1.

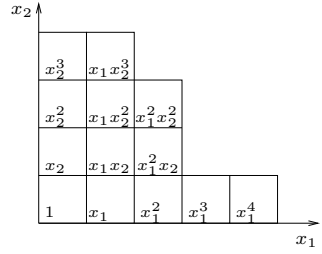
\begin{figure} 
 \center{\input{boxes.pstex_t}} 
\caption{ 
The monomial ideal
$(x_2^4,x_1^2x_2^3,x_1^3x_2,x_1^5)$
determines the partition $(4,4,3,1,1)$.}
\end{figure}

Similarly, the subschemes $S\subset \C^3$ invariant under
the diagonal $\T$-action are in bijective correspondence with
$3$-dimensional partitions. The ideal  
pictured in Figure 2 corresponds to a $3$-dimensional partition
with infinite legs.

In Figures 1 and 2, the
boxes are labelled by the lattice points in the corners
with smallest coordinates,
 a convention which will be followed throughout the paper.

\begin{figure} 
\center{\epsfig{file=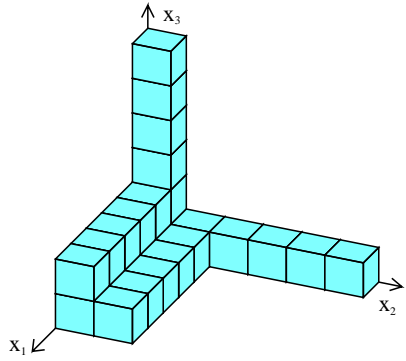, width=95mm}}
\vspace{-43mm}
\caption{The monomial ideal 
$(x_1x_2^2, x_1x_3^2,x_2x_3)$
determines the above $3$-dimensional partition.
The legs in the three coordinate directions are
of infinite length.}
\end{figure}

\subsection{Cohen-Macaulay support}
The first step in the characterization of the restricted
data
\eqref{vvvt}
is to determine the scheme-theoretic support  
$C_\alpha$ of $F_\alpha$. If nonempty, $C_\alpha$ is a
$\T$-invariant, Cohen-Macaulay subscheme of pure dimension 1. 

Let $C\subset \C^3$ be a $T$-fixed 
subscheme of pure dimension 1.
The subscheme $C$ is defined by a 
monomial ideal $$\I_C \subset \C[x_1,x_2,x_3].$$
associated to the 3-dimensional partition $\pi$.
The localisations
$$(\I_C)_{x_1} \subset \C[x_1,x_2,x_3]_{x_1},$$  
$$(\I_C)_{x_2} \subset \C[x_1,x_2,x_3]_{x_2},$$ 
$$(\I_C)_{x_3} \subset \C[x_1,x_2,x_3]_{x_3},$$
are all $T$-fixed, and each corresponds to a 2-dimensional partition 
$\mu^i$.
Alternatively, the 2-dimensional partitions $\mu^i$ can be defined
as the infinite limits of 
the $x_i$-constant cross-sections of $\pi$.
In order for $C$ to have dimension 1,
not all the $\mu^i$ can be empty.

Given a triple $\mmuu=(\mu^1,\mu^2,\mu^3)$ of outgoing partitions,
there exists
a unique {\em minimal} $T$-fixed subscheme 
$$C_{\mmuu}\subset \C^3$$
with outgoing partitions $\mu^i$.
The $3$-dimensional partition corresponding to $C_\mmuu$ is
obtained by taking the union of the infinite cylinders on
the three axes determined by the $2$-dimensional
partitions $\mu^i$.
Let 
\begin{eqnarray*}
\I_{\mu^1}= \mu^1[x_2,x_3] \cdot \C[x_1,x_2,x_3],& \ \ & 
C_{\mu^1}= \O_{\C^3}/\I_{\mu^1},\\
\I_{\mu^2}= \mu^2[x_1,x_3] \cdot \C[x_1,x_2,x_3], & & 
C_{\mu^2}= \O_{\C^3}/\I_{\mu^2}, \\
\I_{\mu^3}= \mu^3[x_1,x_2] \cdot \C[x_1,x_2,x_3], & & 
 C_{\mu^3}= \O_{\C^3}/\I_{\mu^3}.
\end{eqnarray*}
Then $C_\mmuu$ is the union $C_{\mu^1} \cup C_{\mu^2} \cup C_{\mu^3}$
with ideal
$$\I_\mmuu = \bigcap_{i=1}^3 \I_{\mu^i}.$$
If the $\mu^i$ are not all empty, then
$C_{\mmuu}$ is easily seen to be the unique Cohen-Macaulay $T$-fixed curve
in $\C^3$ with these outgoing partitions. By convention, let
$C_{\emptyset,\emptyset,\emptyset}\subset \C^3$
denote the
empty scheme.

Consider the kernel/cokernel sequence associated to the
$\T$-fixed restricted data \eqref{vvvt},
\begin{equation}\label{cvrw}
0 \rightarrow \I_{C_\alpha} \rightarrow \O_{U_\alpha} \stackrel{s}
{\rightarrow} F_\alpha \rightarrow Q_\alpha \rightarrow 0.
\end{equation}
We conclude 
$C_\alpha= C_\mmuu$
where the partitions $\mu^i$ are associated to the edges of $\Delta(X)$
incident to the vertex corresponding to $X_\alpha$.

\label{cmmm}
\subsection{Module $M$}
Since the
support of the quotient $Q_\alpha$ in \eqref{cvrw} is both 0-dimensional
by stability and $\T$-fixed, 
$Q_\alpha$ 
must be supported at the origin.
By Proposition \ref{descl}, the pair $(F_\alpha,s_\alpha)$
corresponds to a $\T$-invariant subsheaf of
$$\lim\limits_{\To}\hom(\m^r,\O_{C_\alpha})/\O_{C_\alpha}
,$$
where 
$\m$  is the
ideal sheaf of the origin in $C_\alpha\subset\C^3$.

Following the notation of Section \ref{cmmm}, let $C_\alpha=C_\mmuu$.
Let
$$M_i = (\O_{C_{\mu^i}})_{x_i}$$
be $\C[x_1,x_2,x_3]$-module  obtained by localisation. So, for instance,
$$
M_1=\C[x_1,x_1^{-1}]\otimes\frac{\C[x_2,x_3]}{\mu^1[x_2,x_3]}\,.
$$
By elementary algebraic arguments,
\begin{eqnarray*}
\lim\limits_{\To}\hom(\m^r,\O_{C_\mmuu}) & 
\cong & \bigoplus_{i=1}^3\lim \limits_{\To}\hom(\m^r,\O_{C_{\mu^i}}) \\
& \cong & \bigoplus_{i=1}^3 M_i.
\end{eqnarray*}
The $\T$-equivariant
$\C[x_1,x_2,x_3]$-module $M_i$ has a canonical $\T$-invariant
element 1. Let
$$
M=\bigoplus_{i=1}^3 M_i.
$$
By Proposition \ref{descl}, the $\T$-fixed pair $(F_\alpha,s_\alpha)$
corresponds to a finitely generated $\T$-invariant
 $\C[x_1,x_2,x_3]$-submodule
\beq{datum}
Q_\alpha\subset M/\langle(1,1,1)\rangle.
\eeq

Conversely, {\em every} finitely generated{\footnote{Here, finitely generated
is equivalent to finite dimensional or Artinian.}}
$\T$-invariant $\C[x_1,x_2,x_3]$-sub\-module
$$Q \subset M/\langle(1,1,1)\rangle$$
occurs as the restriction to $U_\alpha$ of
a $\T$-fixed stable pair on $X$.

\subsection{Box configurations}
We now describe the finitely generated $\T$-invar\-iant
$\C[x_1,x_2,x_3]$-submodules
\beq{datz}
Q \subset M/\langle(1,1,1)\rangle.
\eeq
via labelled box configurations
in the weight space $\mathbb{Z}^3$ of $\T$.

For each of the three partitions $\mu^i$, the module
$M_i$ may be viewed in the space of $\T$-weights as
an infinite cylinder 
$$\text{Cyl}_i\subset \mathbb{Z}^3$$ along the $x_i$-axis with cross
section $\mu^i$. The cylinder extends
in both the positive and negative weight directions.

\begin{figure} 
\center{\epsfig{file=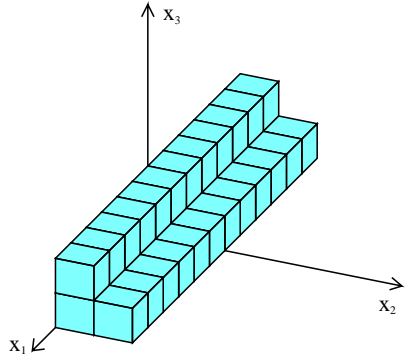, width=95mm}}
\vspace{-37mm}
\caption{
Cylinder associated to the partition $(2,1)$ along the $x_1$-axis
in both the positive and negative directions.}
\end{figure}

The module $M$ is obtained by summing the $M_i$.
For every weight $w$, let ${\mathbf{1}}_w$, ${\mathbf{2}}_w$ and
${\mathbf{3}}_w$ be three independent vectors.
A $\C$-basis for $M$ is
determined by
the set
$$\{ \ {\mathbf{i}}_w\  | \ w\in \text{Cyl}_i \ \}.$$
The $\C[x_1,x_2,x_3]$-module structure on $M$ is clear:
$$x_1 \cdot \mathbf{i}_w = \mathbf{i}_{w+(1,0,0)}, \ \
x_2 \cdot \mathbf{i}_w = \mathbf{i}_{w+(0,1,0)},\ \
x_3 \cdot \mathbf{i}_w = \mathbf{i}_{w+(0,0,1)}.
$$

The union of the cylinders $\text{Cyl}_i$
can be separated into 4 types of weights
$$\bigcup_{i=1}^3 \text{Cyl}_i = 
\mathrm{I}^+\cup\II\cup\III\cup\mathrm{I}^-\subset \mathbb{Z}^3,
$$
where
\begin{enumerate}
\item[$\bullet$]
 $\mathrm{I}^+$ consists of all weights which
have only non-negative coordinates {\em and} which lie in
exactly 1 of the cylinders,
\item[$\bullet$] $\II$ and $\III$ consist of all weights which
lie in exactly 2 and 3 cylinders respectively,
\item[$\bullet$]
 $\mathrm{I}^-$ consists of all weights with at least 1
negative coordinate.
\end{enumerate}

\begin{figure} 
\center{\epsfig{file=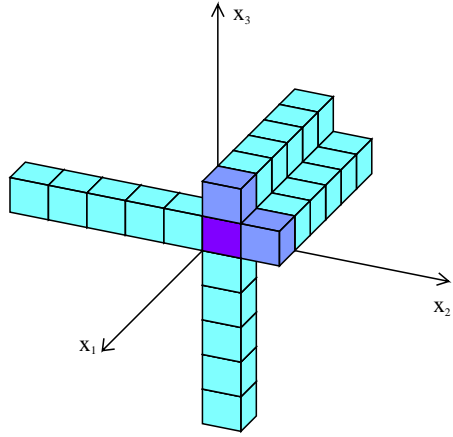, width=95mm}}
\vspace{-37mm}
\caption{Diagram of boxes of type $\mathrm{I}^-$,
$\II$, and $\III$ distinguished by shading
for the partitions $$\hspace{-20pt} \mu^1=(2,1),\ \mu^2=(1),\ \mu^3=(1).$$}
\end{figure}

The submodule $\O_{C_\mmuu}\subset M$ generated by 
$$(1,1,1)=\mathbf{1}_0 + \mathbf{2}_0 + \mathbf{3}_0$$ lies entirely 
in the weight space
$\mathrm{I}^+\cup\II\cup\III$.
The quotient $M/\O_{C_\mmuu}$, described as a $T$-module, is
supported on $\II\cup\III\cup \mathrm{I}^-$ and has the
following $\C$-basis:
\begin{enumerate}
\item[$\bullet$]
If $w\in \mathrm{I}^-$ is
supported on $\text{Cyl}_i$, then 
$$
\C\cdot \mathbf{i}_w \subset M/\O_{C_\mmuu}\,.
$$
\item[$\bullet$]
If $w\in\II$ is supported on  $\text{Cyl}_i$ and $\text{Cyl}_j$, then
$$
\frac{\C\cdot \mathbf{i}_w \oplus\C\cdot
\mathbf{j}_w}{\C\cdot(\mathbf{i}_w+\mathbf{j}_w)}\ \cong\,\C \subset
M/\O_{C_\mmuu}\, .
$$
\item[$\bullet$]
If $w\in\III$, then 
$$
\frac{\C\cdot
\mathbf{1}_w\oplus\C\cdot
\mathbf{2}_w
\oplus\C\cdot
 \mathbf{3}_w}{\C\cdot(1,1,1)_w}\ \cong\,\C^2
\subset
M/\O_{C_\mmuu}\,.
$$
\end{enumerate}
Here, $(1,1,1)_w=\mathbf{1}_w + \mathbf{2}_w + \mathbf{3}_w$.

A finitely generated $T$-invariant $\C[x_1,x_2,x_3]$-submodule
 $$Q\subset M/\O_{C_\mmuu}$$
yields the following  {\em labelled box} configuration
 in $\II\cup\III\cup \mathrm{I}^-$:
a finite number of boxes supported on 
$\II\cup\III\cup \mathrm{I}^-$ where the type $\III$ boxes $w$ 
{\em may} be labelled by an element of
$$
\PP^1=\PP\left(\frac{
\C\cdot \mathbf{1}_w
\oplus\C\cdot \mathbf{2}_w \oplus\C\cdot \mathbf{3}_w}{\C\cdot
(1,1,1)_w}
\right). \vspace{-1mm}
$$
A box signifies the occurance of
the corresponding $\T$-weight in $Q$.
An unlabelled type $\III$ box signifies the
inclusion of the entire 2-dimensional space
$$
\frac{\C\cdot
\mathbf{1}_w\oplus\C\cdot
\mathbf{2}_w
\oplus\C\cdot
 \mathbf{3}_w}{\C\cdot(1,1,1)_w} \subset Q.$$
A labelled type $\III$ box signifies
the inclusion of only the
corresponding 1-dimensional space in
$Q$.

Conversely, given a labelled box configuration,
the following rules ensure that
the corresponding $\T$-submodule $$Q\subset M/\O_C$$ 
is actually a $\C[x_1,x_2,x_3]$-submodule:

\begin{itemize}
\item[(i)] If $w=(w_1,w_2,w_3)\in \mathrm{I}^-$ and if any of 
$$(w_1-1,w_2,w_3),\ (w_1,w_2-1,w_3), \ (w_1,w_2,w_3-1)$$
 support a box then $w$ must support a box.
\item[(ii)] If $w\in\II$,
$w\notin 
\text{Cyl}_i$,  and if any of 
$$(w_1-1,w_2,w_3),\ (w_1,w_2-1,w_3), \ (w_1,w_2,w_3-1)$$
support a box \emph{other than a type $\III$ box labelled by the 
1-dimen\-sional subspace $\C\cdot\mathbf{i}$}, 
then $w$ must support a box.
\item[(iii)] If $w\in\III$ and the span
of the subspaces of 
$$\frac{\C\cdot \mathbf{1}_w\oplus\C\cdot
\mathbf{2}_w\oplus\C\cdot \mathbf{3}_w} 
{\C\cdot (1,1,1)_w}$$ 
induced
by boxes supported on
$$(w_1-1,w_2,w_3),\ (w_1,w_2-1,w_3), \ (w_1,w_2,w_3-1)$$
is nonzero, then $w$ must support a box.
If the span has dimension 1, then $w$ may either support
a box labelled by the span or an unlabelled box.
If the span has dimension 2, then $w$ must support an
unlabelled box.
\end{itemize}
\medskip

The {\em length} of a labelled box configuration is calculated 
by summing the following contributions over the boxes of the
configuration.
Boxes of type $\mathrm{I}^-$ and $\II$ contribute length
1 each. A labelled box of type $\III$ contributes 1 and
an unlabelled box of type $\III$ contributes 2.

A labelled box configuration for $M/O_{C_\mmuu}$ 
is said to have {\em outgoing partitions}
 $\mu^1$, $\mu^2$, and $\mu^3$.

\begin{prop} $\T$-invariant $\C[x_1,x_2,x_3]$-submodules
of $M/\O_{C_\mmuu}$ of length $l$ are in bijective correspondence with
labelled box configurations satisfying (i)-(iii)
with outgoing partitions $\mmuu$ and length $l$. \label{cgh}
\end{prop}

\begin{proof}
Certainly 
$\C[x_1,x_2,x_3]$-submodules satisfy (i)-(iii). An elementary
analysis shows the converse.
\end{proof}

\begin{prop} The {\em reduced} connected components of the  \label{hhh}
moduli space of $\T$-invariant $\C[x_1,x_2,x_3]$-submodules
of  $M/\O_{C_\mmuu}$ are 
products of $\PP^1$s.
\end{prop}

\begin{proof}
Given a labelled box configuration,
let
$L$
be the set labelled type $\III$ boxes. 
Define a {\em path of labelled boxes} to be a sequence of translations
of the form $x_i^{\pm1}$ that stay within $L$.
The set $L$ is divided into disjoint 
path connected subsets.

Two labelled type $\III$ boxes differing by a move of $x_i^{\pm1}$ 
must carry the {\em same} label in $\PP^1$ by rule (iii) above.
Hence, all labelled boxes in each path component of $L$ carry the
same label.

A path component $P\subset L$ is \emph{restricted} if either
of the following two possibilities hold:
\begin{enumerate}
\item[$(\mathsf{r}_+)$]
There is a box in $P$ which is taken
by multiplication by $x_i$
 to a type $\II$ 
box of $M/\O_{C_\mmuu}$ {\em not} occurring in the labelled box
configuration. 
\item[$(\mathsf{r}_-)$] There is a box in $P$ which is taken 
by multiplication by $x_i^{-1}$ to a type $\mathrm{I}^-$
box of the configuration.
\end{enumerate}
The label of such a path component $P$
is forced to be a single point of $\PP^1$ by the rule (iii). 
In the first case above, if the empty type $\II$ box is not in $\text{Cyl}_1$, 
then the label is forced to be $(1,0,0)$.
In the second case, if the type $\mathrm{I}^{-}$ box is in $\text{Cyl}_1$,
then the label is forced to be $(1,0,0)$.

The labellings are the only continuous parameters of the
labelled box configurations.
%$\T$-invariant $\C[x_1,x_2,x_3]$-submodules of  $M/\O_{C_\mmuu}$.
For each unrestricted path component of $L$, 
the label can take any value in $\PP^1$.
Therefore, the moduli space, as a reduced variety, is
simply a product of $\PP^1$s.
\end{proof}

We will use the calligraphic symbol
 $\mathcal{Q}_\mmuu$ to denote components of the moduli space of
$\T$-invariant $\C[x_1,x_2,x_3]$-submodules
of  $M/\O_{C_\mmuu}$. By Propositions \ref{cgh} and \ref{hhh}, the 
components $\mathcal{Q}_\alpha$
correspond to the discrete data of a labelled box configuration ---
forgetting the labelling of the labelled type $\III$ boxes.

\subsection{Local to global}
We have determined the $\T$-fixed restricted data 
$$\O_{U_\alpha} \stackrel{s}{\rightarrow} F_\alpha$$
locally on every $\T$-invariant affine chart $U_\alpha \subset X$.
The gluing condition for different charts is simply the
matching of edge partitions.

We conclude the $\T$-fixed points of  $P_n(X,\beta)$ exactly
arise by distributing
labelled box configurations to
the vertices of $\Delta(X)$,
$$[X_\alpha] \in \mathsf{V}(X) \mapsto Q_\alpha,$$
and partitions to the edges,
$$[C_{\alpha\beta}]\in \mathsf{E}(X) \mapsto \mu_{\alpha\beta},$$
compatible with the outgoing partitions at the
vertices. The vertex data $Q_\alpha$ determines the
edge partitions.

All the moduli in the $\T$-fixed points of $P_n(X,\beta)$
occur at the vertices.
Proposition \ref{hhh} is half of the proof of Theorem \ref{ooo}.
We will complete the proof of Theorem \ref{ooo} by a Zariski
tangent space analysis to show the moduli spaces of
$\T$-invariant  $\C[x_1,x_2,x_3]$-submodules of  $M/\O_{C_\mmuu}$
are scheme-theoretically reduced (and hence nonsingular).

\section{Tangent spaces} \label{tang}
\subsection{$\T$-fixed deformation theory}
The scheme structure on $P_n(X,\beta)$ obtained from
the moduli of stable pairs coincides with the scheme structure
obtained from the moduli of complexes in $D^b(X)$. 
The Zariski tangent space to $P_n(X,\beta)$ at the stable pair
$$I\udot= \left\{ \O_X \stackrel{s}{\rightarrow} F\right\}$$
is $\Ext^0(I\udot,F)$. Derived category $\Ext^0$ may also be
written as $\Hom$.

On each affine chart $U_\alpha \subset X$, the Zariski tangent space to the 
restricted data
$$I\udot_\alpha= \left\{ \O_{U_\alpha} \stackrel{s_\alpha}{\rightarrow} 
F_\alpha\right\}$$
is $\Ext^0(I_\alpha\udot,F_\alpha)$.
There is a global to local restriction map
$$\Ext^0(I \udot,F) \rightarrow 
\bigoplus_\alpha \Ext^0(I_\alpha\udot,F_\alpha)$$  
which need not be an isomorphism. However, if
the stable pair $(F,s)$ is $\T$-fixed, we will
see the induced map
\begin{equation}\label{zwm}
\Ext^0(I \udot,F)^\T \rightarrow 
\bigoplus_\alpha \Ext^0(I_\alpha\udot,F_\alpha)^\T
\end{equation}
 {\em is} an isomorphism.
Here, the superscript $\T$ denotes the $\T$-fixed part, or equivalently,
 the $\T$-weight 0 part.

To complete the proof of 
Theorem \ref{ooo}, we show the Zariski tangent
space $\Ext^0(I_\alpha\udot,F_\alpha)^\T$ to the $\T$-fixed data,
 $$I\udot_\alpha= \left\{ \O_{U_\alpha} \stackrel{s_\alpha}{\rightarrow} 
F_\alpha\right\},$$
described by a labelled
box configuration has dimension equal to the number of unrestricted
path components
of $L_\alpha$, the set of labelled type $\III$ boxes.

The kernel/cokernel sequence \eqref{cvrw} yields the following sequence
of $\T$-modules,
\beq{fixdef}
0\to\Ext^1(Q_\alpha,F_\alpha)\to\Ext^0(I_\alpha\udot,F_\alpha)
\to\Hom(\I_{C_\alpha},F_\alpha).
\eeq
The last term has no 0-weight  piece: the $\T$-weights of
$\I_{C_\alpha}$ lie in the complement
in the weight space  $\Z_{\ge0}^3$ of the $\T$-weights of $F_\alpha$. 
As a result, 
$$\Ext^0(I_\alpha\udot,F_\alpha)^\T
\cong\Ext^1(Q_\alpha,F_\alpha)^\T.$$
Also, the vanishing of the last term shows the $\T$-weight 0 deformations
of the restricted data are supported entirely at the origin. As the
latter can easily be glued,
\eqref{zwm} is an isomorphism.

To avoid calculating with quasi-coherent sheaves 
(or non-finitely generated modules) and passing 
direct limits through derived functors, we work with
a sufficient approximation of $M$,
$$M_i^r=\hom(\m^r,\O_{C_{\mu^i}}), \ \ \ M^r=\oplus_{i=1}^3
M_i^r,$$ 
for $r\gg0$. 
Since $F_\alpha$ is a subsheaf of $M^r$, we obtain
\begin{multline}
\label{fix}
\Hom(Q_\alpha,M^r)\to\Hom(Q_\alpha,M^r/F_\alpha)
\\ \to\Ext^1(Q_\alpha,F_\alpha)\to\Ext^1(Q_\alpha,M^r).
\end{multline}

\begin{lem} \label{fixed}
$\Hom(Q_\alpha,M^r)^\T=\Ext^1(Q_\alpha,M^r)^\T=0$ for $r\gg0$.
\end{lem}

\begin{proof}
By symmetry, we need only prove $\Ext^i(Q_\alpha,M^r_1)$ has no trivial
$T$-subrepresentations for $i=0,1$. 
Since $C_{\mu^1}$ is a product in the $x_1$ direction
and Artinian in the $x_2,x_3$-directions, we have 
$$\m^{r+1}\O_{C_{\mu^1}}=(x_1)\m^r\O_{C_{\mu^1}}, \ \ \ r\gg0,$$
where $(x_1)$ is the ideal of
$\O_{C_{\mu^1}}$ generated by $x_1$. 
The ideal $(x_1)$ is invertible and
abstractly isomorphic as a sheaf to $\O_{C_{\mu^1}}$, 
but twisted by the 1-dimensional
representation of $\T$ with character $t_1$ and associated weight $(1,0,0)$.

Therefore, $M^{r+1}_1\cong  t_1^{-1}\otimes M^r_1$ and
\begin{equation}
\label{xxbg}
\Ext^i(Q_\alpha,M^{r+N}_1)\cong  t_1^{-N} \otimes\Ext^i(Q_\alpha,M^r_1).
\end{equation}
Since $\Ext^i(Q_\alpha,M^r_1)$ is a \emph{finite} sum of 1-dimensional 
$\T$-representations (since $Q$ has 0-dimensional support and 
is finite dimensional),
the $\T$-module \eqref{xxbg}
has no trivial subrepresentations for $N\gg0$.
\end{proof}

Putting together Lemma \ref{fixed}, the previous sequences,
and the stabilization
$$\Hom(Q_\alpha,M/F_\alpha)^\T \cong 
\lim\limits_{\To}
\Hom(Q_\alpha,M^r/F_\alpha)^\T,$$
we conclude that 
$$\Ext^0(I_\alpha\udot,F_\alpha)^\T \cong \Hom(Q_\alpha,M/F_\alpha)^\T$$
for $r\gg0$.

\begin{prop}\label{hhhh}
The
dimension of $\Hom(Q_\alpha,M/F_\alpha)^\T$ equals
the number of unrestricted path components of $L_\alpha$.
\end{prop}

\begin{proof}
The $\T$-weights $w$ of $Q_\alpha$ lie 
in $\mathrm{I}^-\cup\II\cup\III$ and are a subset of the weights of $F_\alpha$.
The latter are a subset of the weights of $M$.
We analyse each type in turn.

If $w$ lies in $\mathrm{I}^-$,
 then $w$  appears in $M$ with multiplicity 1
and so does not appear in $M/F_\alpha$. Thus, the $\T$-weights of $Q$
in $\mathrm{I}^-$ do not contribute to $\Hom(Q,M/F_\alpha)^T$.

If $w$ lies in $\II$, then $w$ appears in both $F_\alpha$ and
$M$ with multiplicity $2$. Again,
 $w$ does not appear in $M/F_\alpha$ and so does not contribute to
$\Hom(Q,M/F_\alpha)^T$.
Similarly, if $w$ lies in $\III$ with multiplicity 2, then $w$
 appears in both $F_\alpha$ and $M$ with multiplicity 3
and does not contribute to
$\Hom(Q_\alpha,M/F_\alpha)^T$.

If $w$ lies in $\III$ with multiplicity 1, then $w$ appears in $F_\alpha$ 
with multiplicity 2 but in $M$  with multiplicity 3.
The multiplicity of $w$ in $M/F_\alpha$  is 1. 
Thus, we find an at most a 1-dimensional subspace of 
$\Hom(Q_\alpha,M/F_\alpha)^T$ corresponding to such $w$. 
However, the
analysis used in the proof of 
Proposition  \ref{hhh} shows that the $\C[x_1,x_2,x_3]$-module
structure forces any morphism in the $w$-box 
to be equal to the morphism in
any other box in the same path component of $L_\alpha$. And,
 if the path component is restricted,
then the morphism is 0 over the whole path component.

Therefore, $\Hom(Q_\alpha,M/F_\alpha)^\T$ 
has dimension equal to the number of unrestricted
path components of $L_\alpha$.
\end{proof}

Propositions \ref{hhh} and \ref{hhhh} imply Theorem \ref{ooo}.
Proposition \ref{hhh} provides a description of the {\em reduced}
$\T$-fixed components of $P_n(X,\beta)$. 
The $\T$-fixed Zariski tangent space obtained
from Proposition \ref{hhhh} establishes the nonsingularity
of the scheme structure.

Let $\mathbf{Q}\subset P_n(X,\beta)$ be a component of the
$\T$-fixed locus. We have proven
$$\mathbf{Q} = \prod_{[X_\alpha]\in \bV(X)} \mathcal{Q}_\alpha$$
where $\mathcal{Q}_\alpha$ is a component of moduli space of
labelled box configurations. Each element of $\mathbf{Q}$
can be described by a labelled box configuration
$$[X_\alpha]\in \bV(X) \mapsto Q_\alpha$$
at each vertex.
We will follow the above notation throughout the paper.

\label{xew}
\subsection{$\T_0$-fixed deformation theory} \label{tocy}
Let $X$ be a toric Calabi-Yau 3-fold with canonical form 
$$\omega\in H^0(X,K_X)$$
invariant under a 2-dimensional subtorus
$\T_0\subset \T$.
On the affine chart 
$$\C^3\cong U_\alpha\subset X,$$
 the subtorus $\T_0$
must act trivially on the form $dx_1 \wedge dx_2 \wedge dx_3$
and hence must be defined by
$$\T_0 = \{ \ (t_1,t_2,t_3) \in \T \ | \ t_1t_2t_3=1 \ \}.$$

Consider the $\T$-fixed restricted data studied above,
$$I\udot_\alpha= \left\{ \O_{U_\alpha} \stackrel{s_\alpha}{\rightarrow} 
F_\alpha\right\}.$$
Certainly $I\udot_\alpha$ is also $\T_0$-fixed.

\begin{lem} $\Hom(\I_{C_\alpha},F_\alpha)$ 
contains no $\T_0$-fixed representation. \label{ffqw}
\end{lem}

\begin{proof} 
The space $\Hom(\I_{C_\alpha},F_\alpha)^{\T_0}$ may be decomposed as a direct
sum of weight spaces for the quotient torus
$\C^*\cong \T/\T_0$.
Homomorphisms 
%$$\phi\in \Hom(\I_{C_\alpha},F_\alpha)^{\T_0}$$
of $\C^*$-weight $w\in \Z$
multiply the $\T$-submodules of
$\I_{C_\alpha}$ by $(x_1x_2x_3)^w$ to give $\T$-submodules of $F_\alpha$. 

We have seen in Section \ref{xew} 
there are no such homomorphisms of $\C^*$-weight $0$. 
The same argument shows the non-existence in the positive 
weight case. The $\T$-weights of $(x_1x_2x_3)^{w\geq 0}\I_{C_\alpha}$ 
all lie in the complement in $\Z^3_{\ge0}$
of the $\T$-weights of $F_\alpha$.
%Let $\phi$ be homomorphism of $\C^*$-weight $w<0$.
Since $$F_\alpha\subset M=\bigoplus _{i=1}^3 M_i,$$
the vanishing of all $\T_0$-fixed homomorphisms of $\C^*$-weight $w<0$
between $\I_{C_\alpha}$ and $M_i$ implies the Lemma. 

\begin{figure} \center{\input{moves.pstex_t}}
\caption{The partition $\mu^1[x_2,x_3]$. \label{moves}}
\end{figure}
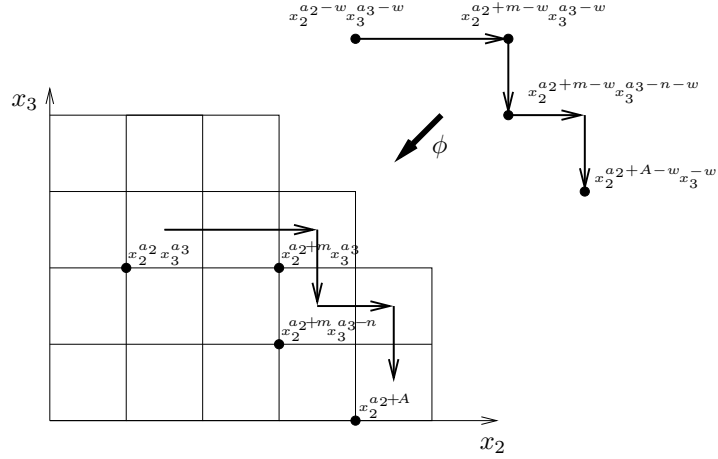

In fact, by
symmetry, we need only study $M_1$.
We write 
$$M_1=\C[x_1,x_1^{-1}]\otimes_\C \frac{\C[x_2,x_3]}{\mu^1[x_2,x_3]}$$ 
for the outgoing partition
$\mu^1$. 
Let $\phi\in \Hom(\I_{C_\alpha},M_1)^{\T_0}
$ be homomorphism of $\C^*$-weight $w<0$.
Pick a nonzero element $$x_1^{a_1}x_2^{a_2}x_3^{a_3} 
\in \text{Im}(\phi) \subset M_1.$$ 
Then, perhaps after scaling, we have 
$$\phi(x_1^{a_1-w}x_2^{a_2-w}x_3^{a_3-w})= x_1^{a_1}x_2^{a_2}x_3^{a_3},$$ where
$x_1^{a_1-w}x_2^{a_2-w}x_3^{a_3-w}
\in \I_{C_\alpha}$.

By multiplying by $x_1^N$ if necessary, we may assume that $a_1\gg0$. 
In the limit, we may see the problem as essentially
two dimensional in the variables $x_2,x_3$.

There is a maximal $m\ge0$ satisfying 
$$0\neq x_2^{a_2+m}x_3^{a_3}\in \frac{\C[x_2,x_3]}{\mu^1[x_2,x_3]}\,.$$
Such a monomial $x_2^{a_2+m}x_3^{a_3}$ is said to be
{\em $x_2$-maximal}. Consider the set 
$$S = \{ \ n'\ge0\ | \ x_2^{a_2+m}x_3^{a_3-n'} {\text{ \ is not 
$x_2$-maximal}}\ \}.$$
Let $n$ be the minimal element of $S$. If $S$ is empty, let 
$n=a_3$. Unless $a_3=0$, $n$ must be positive.
See Figure \ref{moves}.

By the minimality of $n$ and the strict negativity of $w$, 
$$x_2^{a_2+m-w}x_3^{a_3-n-w} \in \mu^1[x_2,x_3] \subset \C[x_2,x_3].$$
Since 
$a_1\gg0$, we also have 
$x_1^{a_1-w}x_2^{a_2+m-w}x_3^{a_3-n-w}\in\I_{C_\alpha}$.
Then,
 $$x_3^n\phi(x_1^{a_1-w}x_2^{a_2+m-w}x_3^{a_3-n-w})=x_1^{a_1}x_2^{a_2+m}
x_3^{a_3}$$
is nonzero in $M_1$. 
We find
$$\phi(x_1^{a_1-w}x_2^{a_2+m-w}x_3^{a_3-n-w})=
x_1^{a_1}x_2^{a_2+m}x_3^{a_3-n}\neq 0.$$
We have found another nonzero
element of $M_1$ in the image of $\phi$ with smaller $x_3$ exponent.

Inductively, we reduce the $x_3$ exponent to 0, see
Figure \ref{moves}. 
We find there is an $A\ge0$ for which
$x_1^{a_1-w}x_2^{a_2+A-w}x_3^{-w}\in\I_C$ and
$$\phi(x_1^{a_1-w}x_2^{a_2+A-w}x_3^{-w})=x_1^{a_1}x_2^{a_2+A}$$
 is nonzero in
$M_1$. 
After multiplying by a non-negative  power of $x_2$,
we may also assume $x_2^{a_2+A}$ to be  $x_2$-maximal. 
Hence,
$$x_2^{a_2+A-w}\in\mu^1[x_2,x_3], \ \  \
 x_1^{a_1-w}x_2^{a_2+A-w}\in \I_C.$$ 
We now obtain a contradiction since
$$
\phi(x_1^{a_1-w}x_2^{a_2+A-w}x_3^{-w})
=x_1^{a_1}x_2^{a_2+A}=x_3^{-w}\phi(x_1^{a_1-w}x_2^{a_2+A-w})
$$
and $\phi(x_1^{a_1-w}x_2^{a_2+A-w})=0$ since its
third $\T$-weight is negative.
\end{proof}

Lemma \ref{ffqw} and the $\T_0$-fixed part of
sequence \eqref{fixdef}
imply the local to global sequence
is an isomorphism for $\T$-fixed stable pairs.
\begin{equation}
\label{vvhl}
\Ext^0(I \udot,F)^{\T_0} \cong
\bigoplus_\alpha \Ext^0(I_\alpha\udot,F_\alpha)^{\T_0}
\end{equation}
Also, since the proof of 
Lemma \ref{fixed} is valid for $\T_0$ in place of $\T$, we obtain
the following result.

\begin{lem}\label{gg23}
The Zariski tangent space $\Ext^0(I_\alpha\udot,F_\alpha)^{\T_0}$ 
to the $\T_0$-fixed
locus 
is equal to $\Hom(Q_\alpha,M/F_\alpha)^{\T_0}$.
\end{lem}

\subsection{Nonsingularity} \label{ccooj1}
The question of the nonsingularity of the $\T_0$-fixed loci
$$P_n(X,\beta)^{\T_0} \subset P_n(X,\beta)$$
is very natural and plays a crucial role in our 
study of the Calabi-Yau vertex.

\begin{conj} \label{vvt}
The loci $P_n(X,\beta)^{\T_0}$ are nonsingular.
\end{conj}

By the local to global relation of tangent spaces \eqref{vvhl},
Conjecture \ref{vvt} is equivalent to the nonsingularity of the
moduli space of 
$\T_0$-invariant submodules of $M/O_{C_\mmuu}$.

A {\em local toric surface} is the total space of a toric line
bundle over a toric surface. If $X$ is a local toric surface,  
The restricted supports $C_\alpha$
have only 
1 or 2 legs. By Propositions \ref{hhh} and \ref{hhhh}, the
loci $P_n(X,\beta)^\T$ are isolated points. 
Lemma \ref{gg23} easily implies the
tangent space $\Ext^0(I_\alpha\udot,F_\alpha)^{\T_0}$
to the restricted data is 0, so the $\T_0$-fixed
Zariski tangent space is no larger.
Hence Conjecture \ref{vvt} is proven in the local toric surface
case.

Conjecture \ref{vvt} would follow in general if
the tangent vectors 
$$\Hom(Q_\alpha,M/F_\alpha)^{\T_0}$$
to the restricted data could be exponentiated. 
The latter is an essentially combinatorial condition which holds
in all the examples we have studied.

\subsection{Dimension}
\label{ccooj2}

For any finite $\T$-module $K$, let 
$$\chi\_K(t_1,t_2,t_3)\in \Z[t_1^{\pm},t_2^\pm,t_3^\pm]$$ be the
$\T$-character,
and let
$$\overline{\chi}_K(t_1,t_2,t_3) = \chi\_K(t_1^{-1}, t_2^{-1}, t_3^{-1}).$$
Define the Laurent polynomial $\gamma\_K$ by
$$\gamma\_K= \chi\_K -
\chi\_K \bar{\chi}_K \frac{(1-t_1)(1-t_2)}{t_1t_2}\,.$$
Finally, let
$$\gamma\_{K,0}(t_1,t_2) = \gamma\_K(t_1,t_2, t_1^{-1}t_2^{-1})\in
\C[t_1^\pm, t_2^\pm].$$
Let $\mathsf{Con}(K,0)$ be the {\em constant} term of $\gamma\_{K,0}$.

Given $\T$-fixed restricted data $(F_\alpha, s_\alpha)$,
let $F_\alpha^{c}$ be the finite length
$\T$-module obtained
by cutting off{\footnote{The cut-offs may be
taken to be simple cuts perpendicular to the
axes. Since we are only interested
in suitably large cut-offs, cuts which are finitely
jagged are also fine to consider as well.}}
the infinite legs of $F_\alpha$ in the
3 positive directions. 
The parity of the dimension of
the $\T_0$-fixed
tangent space is determined by the following result.

\begin{conj} \label{vvgg}
For all suitably large cut-offs $F^c_\alpha$,
$$\text{\em dim}_\C \ \Ext^0(I_\alpha \udot,F_\alpha)^{\T_0} =
\mathsf{Con}(F^{c}_\alpha,0) \ \mod 2.$$
\end{conj}

If the support curve $C_\alpha$ has 1 or 2 legs, 
we have seen the dimension of 
  $\Ext^0(I_\alpha \udot,F_\alpha)^{\T_0}$ is 0. 
Conjecture \ref{vvgg} then asserts that
$\mathsf{Con}(F^{c}_\alpha,0)$ is {\em even}.
The proof is easily obtained by box removal in the manner
described in Section 4.11 of \cite{MNOP1}.
We leave the details to the reader.

If the  parity of the dimensions of the Zariski tangent
spaces to $P_n(X,\beta)$ were {constant} in families of stable
pairs, 
Conjecture \ref{vvgg}
would be very natural. It is tempting to look for a formal proof
of the 
even jumping of tangent dimensions via an appropriate 
 symplectic form in the ambient geometry. 

Unfortunately, we know
the conservation of the parity of the Zariski tangent spaces is
{\em false} in general. An example can be found 
in a $K3$-fibration,
$$\epsilon:X\rightarrow \PP^1,$$
where the fiber $X_0$ over $0\in \PP^1$
is
a $K3$ surface with a fully obstructed{\footnote{Here, we
require the algebraic classes $[L_1]$, $[L_2]$, and
$[L_1]+[L_2]$ are all obstructed to first order in the family $\epsilon$.}
$A_2$-configuration,
$$L_1 \cup L_2 \subset X_0.$$
Here, $L_1$ and $L_1$ are $(-2)$-curves meeting transversely in a
point,
$$\chi( \O_{L_1\cup L_2}) =1.$$
The Zariski tangent spaces of $P_2(X,[L_1]+[L_2])$ are 
seen to jump from 1-dimension to 2-dimensions.{\footnote{
A similar odd jumping phenomena can be found
in the moduli space \\ $I_2(X,[L_1]+[L_2])$ of
ideal sheaves.}}

The sign of Conjecture \ref{vvgg} is therefore not purely
formal and depends essentially on toric geometry. In fact, the
sign may be viewed as the only aspect of the
Calabi-Yau vertex calculation which is not purely formal.{\footnote{A
position taken by J. Bryan.}}

\section{Equivariant vertex} \label{eqver}

\subsection{Localization} \label{vl}
Let $X$ be a nonsingular, quasi-projective, toric $3$-fold.
The $\T$-equivariant obstruction theory of $P_n(X,\beta)$ obtained
from deformations in $D^b(X)$ admits a 
2-term $\T$-equivariant free resolution
\begin{equation}\label{oobb}
E^{-1} \rightarrow E^0,
\end{equation}
see \cite{pt}.
Denote the dual of the restriction of \eqref{oobb} to a
component $\mathbf{Q}$ of 
the $\T$-fixed locus 
\begin{equation*} 
\iota:\mathbf{Q}\hookrightarrow P_n(X,\beta),
\end{equation*}
by
$E_{\mathbf{Q},0} \rightarrow E_{\mathbf{Q},1}$.

The $\T$-action on \eqref{oobb} restricts to a fiberwise
$\T$-action. Let
$E_{\mathbf{Q},i}^{\T}$ and $E_{\mathbf{Q},i}^m$ denote the
sub-bundles with 0 and nonzero characters respectively.
A $\T$-fixed obstruction theory 
on $\mathbf{Q}$ is obtained from
$$E_{\mathbf{Q},0}^{\T} \rightarrow E_{\mathbf{Q},1}^{\T}\,.$$
Let $[\mathbf{Q}]^{vir}$ be the associated virtual
class \cite{GraberP}.

By the virtual localization formula
\cite{GraberP} we have
\begin{equation*}
\sum_{\mathbf{Q}\subset P_n(X,\beta)} \iota_*\left( \frac
{e(E^m_{\mathbf{Q},1})}{e(E^m_{\mathbf{Q},0})}
\cap [\mathbf{Q}]^{vir} \right)
= [P_n(X,\beta)]^{vir} \in A^{\mathbf{T}}_*(P_n(X,\beta))_{loc}\,.
\end{equation*}
Here, $e$ denotes the $\T$-equivariant Euler class, or equivalently,
the top Chern class. The formula
 holds in the localized $\T$-equivariant
Chow ring of $P_n(X,\beta)$.

\subsection{$\T$-fixed obstruction theory}\label{fxy}
The data of the $\T$-fixed obstruction theory on $\mathbf{Q}$ 
includes a
morphism $\phi$ to the cotangent complex,
$$[(E_{\mathbf{Q}}^{-1})^\T \rightarrow (E^0_{\mathbf{Q}})^\T] 
\stackrel{\phi}{\rightarrow}
L^\bullet_{\mathbf{Q}}\ ,$$
for which
$h^0(\phi)$ is an isomorphism and
$h^{-1}(\phi)$ is a surjection.
Since $\mathbf{Q}$ is nonsingular by Theorem \ref{ccc}, the 2-term
cut-off of the cotangent complex of $\mathbf{Q}$ can be taken simply
to be
$$L^\bullet_{\mathbf{Q}} = [ 0 \rightarrow \Omega_{\mathbf{Q}}].$$
We can, after exchanging $E^\bullet$, assume the morphism $\phi$ is
represented by a map of complexes. Then, the sequence
$$(E_{\mathbf{Q}}^{-1})^\T \rightarrow (E^0_{\mathbf{Q}})^\T
\stackrel{\phi^0}{\rightarrow} \Omega_{\mathbf{Q}} \rightarrow 0$$
is exact. The kernel on the left is a bundle $K$, and
$$[\mathbf{Q}]^{vir} = e(K^\vee) \cap [\mathbf{Q}]$$
by the definition of the virtual class \cite{bf}.

Let $T_{\mathbf{Q}}$ denote the tangent bundle of $\mathbf{Q}$.
The obstruction bundle $K^\vee$ is determined in $K$-theory by
\begin{equation}
\label{xrby}
[K^\vee]= [E_{\mathbf{Q},1}^\T]-[E_{\mathbf{Q},0}^\T]+[T_{\mathbf{Q}}]
\end{equation}
The formal Euler class expression
$$e(T_{\mathbf{Q}}) \cdot \frac{e(E_{\mathbf{Q},1}^\T)}{
e(E_{\mathbf{Q},0}^\T)}
\in A^*(\mathbf{Q})$$
is therefore well-defined and equal to $e(K^\vee)$.

Together with the localization formula, we obtain the following
result,
\begin{equation}\label{gpf}
\sum_{\mathbf{Q}\subset P_n(X,\beta)} \iota_*\left(e(T_{\mathbf{Q}}) 
\cdot \frac
{e(E_{\mathbf{Q},1})}{e(E_{\mathbf{Q},0})}
\cap [\mathbf{Q}] \right)
= [P_n(X,\beta)]^{vir} 
\end{equation}
in
$A^{\mathbf{T}}_*(P_n(X,\beta))_{loc}$.

\subsection{Local to global}
We will calculate  the difference
\begin{equation}\label{ssdd}
[E_{\mathbf{Q},0}] - [E_{\mathbf{Q},1}]
\end{equation}
in the $\T$-equivariant $K$-theory of $\mathbf{Q}$.
Consider a stable pair
$$I\udot= \left\{ \O_X \stackrel{s}{\rightarrow} F\right\}, \ \ 
[I\udot] \in \mathbf{Q}.$$
The difference \eqref{ssdd} restricted to $[I\udot]\in \mathbf{Q}$
is simply
the virtual tangent 
space,{\footnote{Tracelessness is automatic for toric varieties.}}
\begin{equation}
\label{vvtan}
\mathcal{T}_{\left[{I}\udot\right]} =
\Ext^1(I\udot,I\udot)- \Ext^2(I\udot,I\udot).
\end{equation}
Our first goal is a canonical calculation of the $\T$-representation
\eqref{vvtan} relatively over $\mathbf{Q}$.

By \cite{pt}, the virtual tangent space at $I\udot$
is given by
$$
\mathcal{T}_{\left[{I}\udot\right]} 
%\extt^{1}({I}\udot,{I}\udot) - \extt^{2}({I}\udot,{I}\udot)
= \chi(\mathcal{O}, \mathcal{O}) - \chi({I}\udot, {I}\udot)
$$
where
$$
\chi({F}\udot, {G}\udot) =
\sum_{i=0}^{3}(-1)^{i}\extt^{i}({F}\udot,{G}\udot)\,.
$$
We can compute each
Euler characteristic using the local to global spectral sequence
\begin{align*}
  \chi({I}\udot, {I}\udot)
&= \sum_{i,j=0}^{3}(-1)^{i+j}H^{i}(\sheafext^{j}({I}\udot,
{I}\udot)) \\
&= \sum_{i,j = 0}^{3}(-1)^{i+j}
\mathfrak{C}^{i}(\sheafext^{j}({I}\udot,{I}\udot))
\,,
\end{align*}
where, in the second line,
we have replaced the cohomology terms with the
Cech complex with respect to the open affine cover $\{U_\alpha\}$.
Though the modules of the Cech complex are infinite-dimensional,
they have finite-dimensional weight spaces and, therefore,
their $\T$-character is well defined as a formal power
series.

Since $F$ is supported on
the $\T$-invariant
 curves corresponding to the edges of $\Delta(X)$,
we have $I\udot = \O_X$
on the intersection of three or more $U_{\alpha}$.
Therefore, only the $\mathfrak{C}^{0}$ and 
$\mathfrak{C}^{1}$ terms contribute to
the calculation. We find,
\begin{multline}\label{virtT}
  \mathcal{T}_{\left[{I}\udot\right]} =
\bigoplus_{\alpha}\left(\Gamma(U_{\alpha})
- \sum_{i}(-1)^{i}\Gamma(U_{\alpha}, 
\sheafext^{i}({I}\udot,{I}\udot))\right) \\
-
\bigoplus_{\alpha, \beta}\left(\Gamma(U_{\alpha\beta}) -
\sum_{i}(-1)^{i}\Gamma(U_{\alpha\beta}, 
\sheafext^{i}({I}\udot,{I}\udot))\right)
\,.
\end{multline}

The calculation of $\mathcal{T}_{\left[{I}\udot\right]}$
is reduced to a sum over all the vertices and edges of 
$\Delta(X)$.
We have the restrictions of $I\udot$,
$$I_\alpha \udot \in D^b(U_\alpha), \ \ I_{\alpha\beta}\udot
\in D^b(U_{\alpha\beta}),$$
and we need to compute
\begin{eqnarray*}
\Gamma(U_\alpha)& -& 
\sum_{i}(-1)^{i}\extt^i(I\udot_\alpha,I\udot_\alpha), \\
\Gamma(U_{\alpha\beta})& -& 
\sum_{i}(-1)^{i}\extt^i(I\udot_{\alpha\beta},I\udot_{\alpha\beta}).
\end{eqnarray*}
%In the vertex case, $\Gamma(U_\alpha) \cong \mathbf{C}[x,y,z]$, and
%in the edge case, $\Gamma(U_\alpha \cap U_\beta) \cong
%\mathbf{C}[x,y,z,z^{-1}]$. 
The vertex and edge cases will be treated 
separately.

A similar strategy was pursued in \cite{MNOP1} to calculate
the DT invariants of toric 3-folds. The difference here
occurs entirely in the vertex terms. The
edge terms are identical.

\subsection{Vertex calculation}
\label{vc}
Let $R$ be the coordinate ring,
$$
R = \C[x_1,x_2,x_3] \cong \Gamma(U_{\alpha}).
$$
Following
the conventions of Section \ref{tgeom}, the $\T$-action on $R$ is 
\begin{equation*}
  (t_1,t_2,t_3)\cdot x_i = t_i x_i  \,.
\end{equation*}

Let $\Ib_\alpha\udot$ denote the universal complex on
$\mathcal{Q}_\alpha \times U_\alpha$.
Consider a $\T$-equivariant free
resolution of $\Ib_\alpha\udot$,
\begin{equation}
  \label{resol}
 \{ \F_{s} \rightarrow \dots \rightarrow \F_{1}\}
\cong 
\Ib_\alpha\udot \ \in D^b(\mathcal{Q}_\alpha \times U_\alpha).
\end{equation}
Each term in \eqref{resol}
can be taken to have the form
$$
\F_i = \bigoplus_j \mathcal{L}_{ij}
\otimes R(d_{ij})\,, \quad d_{ij} \in \Z^3,$$
where $\mathcal{L}_{ij} \in \text{Pic}(\mathcal{Q}_\alpha)$.
The Poincar\'e polynomial
$$
P_\alpha = \sum_{i,j} (-1)^i \, [\mathcal{L}_{ij}] \otimes t^{d_{ij}}
\ \in    K(\mathcal{Q}_\alpha)    \otimes_{\mathbb{Z}} 
\Z[t_1^\pm,t_2^\pm,t_3^\pm]$$
does not depend on the choice of the resolution
\eqref{resol}. 
We require the $K$-theoretic data to keep track of twisting over
$\mathcal{Q}_\alpha$.

We denote the $\T$-character of  $F_\alpha$ by $\FFF_\alpha$.
By Theorem \ref{ccc} and the
sequence
$$0 \rightarrow \O_{C_\alpha} \rightarrow F_\alpha \rightarrow Q_\alpha 
\rightarrow 0,$$
we have a complete understanding of
the representation $\FFF_\alpha$. 
The $\T$-eigenspaces of $F_\alpha$ correspond to the
$\T$-eigenspaces of $\O_{C_\alpha}$ and 
the boxes of the labelled configuration
 associated to $Q_\alpha$. 
The boxes contributes
monomials to $\FFF_\alpha$.
For each unrestricted path component $P\subset L_\alpha$,
let $\PP_P$ denote the associated factor of $\mathcal{Q}_\alpha$.
For every labelled box in $P$, 
we tensor the corresponding  character  {by} $\O_{\PP_P}(-1)$.
The $K$-theory factors associated to $\O_{C_\alpha}$ and the
unlabelled boxes are trivial.
Of course, an unlabelled type $\III$ box contributes twice.
The result determines 
$$\FFF_\alpha \in K(\mathcal{Q}_\alpha) 
\otimes_{\mathbb{Z}} \Z(t_1,t_2,t_3).$$
The rational dependence on the $t_i$ is elementary.

From
the resolution \eqref{resol}, we see
that the Poincar\'e polynomial
$P_\alpha$ is related to the $\T$-character
of $F_\alpha$ as follows:
\begin{equation}
\FFF_\alpha =
\frac{1+P_\alpha}{(1-t_1)(1-t_2)(1-t_3)} \label{PQ}
\,.
\end{equation}
Hence, we may effectively compute $P_\alpha$.

The family
$\chi(\Ib_\alpha\udot,\Ib_\alpha\udot)$ of virtual
representations  
over $\mathcal{Q}_\alpha$ is given by the
following alternating sum
\begin{align*}
\chi(\Ib_\alpha\udot,\Ib_\alpha\udot)  &= \sum_{i,j,k,l} (-1)^{i+k}
\mathcal{L}_{ij} \otimes \mathcal{L}_{kl}^\vee
\otimes \Hom_R(R(d_{ij}), R(d_{kl}))
\\
&=  \sum_{i,j,k,l} (-1)^{i+k}\mathcal{L}_{ij} \otimes \mathcal{L}_{kl}^\vee
\otimes
R(d_{kl}-d_{ij})\,.
\end{align*}
Therefore, the $\T$-character is 
$$
\tr_{\chi(\Ib_\alpha,\Ib_\alpha)} =
\frac{P_\alpha  \,\overline{P}_\alpha}
{(1-t_1)(1-t_2)(1-t_3)} \,.
$$
The dual bar operation 
$$\gamma \in K(\mathcal{Q}_\alpha) \otimes_\Z \Z(\!(t_1,t_2,t_3)\!) \mapsto
\overline{\gamma} \in K(\mathcal{Q}_\alpha) \otimes_\Z 
\Z(\!(t_1^{-1},t_2^{-1},
t_3^{-1})\!)$$
is negation on $K(\mathcal{Q}_\alpha)$ and
$$t_i \mapsto t_i^{-1}$$
on the variables.

We find the $\T$-character of 
the $\alpha$ summand of  virtual tangent
space $\mathcal{T}_{\left[{I}\udot\right]}$
in  \eqref{virtT} is
$$
\frac{1-P_\alpha \, \overline{P}_\alpha}
{(1-t_1)(1-t_2)(1-t_3)} \,.
$$
Using \eqref{PQ}, we may express the answer in terms of
$\FFF_\alpha$,
\begin{equation}\label{vertexchar}
  \tr_{R-\chi(\Ib\udot_\alpha,\Ib_\alpha\udot)} 
= \FFF_{\alpha} -
\frac{\overline{\FFF}_\alpha}{t_1t_2t_3} +  \FFF_{\alpha}
\overline{\FFF}_\alpha \frac{(1-t_1)(1-t_2)(1-t_3)}{t_1 t_2 t_3} \,.
\end{equation}
On the right side of
 \eqref{vertexchar}, the rational functions
should be expanded 
in ascending powers in the $t_i$.

\subsection{Edge calculation}
We now consider the summand of \eqref{virtT} corresponding to
a pair $(\alpha,\beta)$. Our calculations will involve modules over the ring
$$
R=\Gamma(U_{\alpha\beta})  =
 \C[x_2,x_3]\otimes_\C
\C[x_1,x_1^{-1}] \,.
$$
The $\C[x_1,x_1^{-1}]$ factor will result
only in the overall factor
$$
\delta(t_1) = \sum_{k\in \Z} t_1^k,
$$
the formal $\delta$-function at $t_1=1$, in the $\T$-character.
Let
$$
\FFF_{\alpha\beta}  = \sum_{(k_2,k_3) \in \mu_{\alpha\beta}} t_2^{k_2}
t_3^{k_3}
$$
be the generating function for the edge partition $\mu_{\alpha\beta}$.
Arguing as in the vertex case, we find
\begin{multline}\label{edgechar}
 - \tr_{R-\chi(I_{\alpha\beta},I_{\alpha\beta})} \\
= \delta(t_1) \left(- \FFF_{\alpha\beta} -
\frac{\overline{\FFF}_{\alpha\beta}}{t_2 t_3} +  \FFF_{\alpha\beta}
\overline{\FFF}_{\alpha\beta} \frac{(1-t_2)(1-t_3)}{t_2 t_3} \right)\,.
\end{multline}
Because of the relations
$$
\delta(1/t) = \delta(t) = t\delta(t) \,,
$$
the character \eqref{edgechar} is invariant under the change
of variables \eqref{trans}.

Since there is no $\T$-fixed moduli away from the vertices, there
is no $K$-theoretic data associated to the edges.

\subsection{Redistribution}
\label{seqv}
We
now  redistribute the terms of the vertex \eqref{vertexchar}
and edge \eqref{edgechar}
contributions  so both become 
Laurent {\em polynomials} in the variables $t_i$.

The edge character \eqref{edgechar} can be written as
\begin{equation}
  \label{FF}
  \frac{\GGG_{\alpha\beta}(t_2,t_3)}{1-t_1} +
t_1^{-1}\frac{\GGG_{\alpha\beta}(t_2,t_3)}{1-t_1^{-1}} \,,
\end{equation}
where the first 
term in \eqref{FF} is expanded
in ascending powers of $t_1$, and the second term is
expanded in descending powers. Here
$$
\GGG_{\alpha\beta} = - \FFF_{\alpha\beta} -
\frac{\overline{\FFF}_{\alpha\beta}}{t_2 t_3} +  \FFF_{\alpha\beta}
\overline{\FFF}_{\alpha\beta} \frac{(1-t_2)(1-t_3)}{t_2 t_3} \,.
$$

Define a new vertex character $\bV_\alpha$ by the following
modification,
\begin{equation}\label{gx34}
\bV_\alpha = \tr_{R-\chi(\Ib\udot_\alpha,\Ib_\alpha\udot)} 
 + \sum_{i=1}^3
\frac{\GGG_{\alpha\beta_i}(t_{i'},t_{i''})}{1-t_i}\,,
\end{equation}
where $\beta_1, \beta_2,\beta_3$  are the three neighboring
vertices
and $$\{t_i, t_{i'}, t_{i''}\} = \{t_1,t_2,t_3\}.$$
The character $\bV_\alpha$ depends
{\em only on the local data $\mathcal{Q}_\alpha$}.
%The definition of the character $\bV_\beta$ for other fixed points
%$X_\beta \in X^\T$ involves changing the variables as in \eqref{trans}.
Similarly, we define
 \begin{equation*}
   \bE_{\alpha\beta} =
t_1^{-1}\frac{\GGG_{\alpha\beta}(t_2,t_3)}{1-t_1^{-1}}
- \frac{\GGG_{\alpha\beta}\left(t_2\, t_1^{-m_{\alpha\beta}},
t_3\, t_1^{-m'_{\alpha\beta}}\right)}{1-t_1^{-1}}\,.
\end{equation*}
The term $\bE_{\alpha\beta}$ is canonically associated to the edge.
Formulas \eqref{vertexchar} and \eqref{edgechar} yield the following
result.

\begin{thm} The $\T$-character of $\mathcal{T}_{\left[\Ib\udot\right]}$
over $\mathbf{Q}$
is given by
\begin{equation}
  \label{charTvirt}
  \tr_{\mathcal{T}_{\left[\Ib\udot\right]}} =
\sum_{[X_\alpha]\in \bV(X)}\ \bV_\alpha + 
\sum_{[C_{\alpha\beta}]\in \bE(X)} \bE_{\alpha\beta} \,.
\end{equation}
\end{thm}

\begin{lem} Both $\bV_\alpha$ and $ \bE_{\alpha\beta}$
are Laurent polynomials in the $t_i$.
\end{lem}

\begin{proof}
  The numerator of $\bE_{\alpha\beta}$ vanishes at $t_1=1$ and
therefore is divisible by the denominator. The claim for $\bV_\alpha$
follows from
$$
\FFF_\alpha = \sum_{i=1}^3\frac{\FFF_{\alpha\beta_i}}{1-t_1} + \dots\,,
$$
where $\beta_1$, $\beta_2$, and $\beta_3$ are
the neighboring vertices and
the dots stand for terms regular at $t_1=1$.
\end{proof}

\subsection{The equivariant vertex}
Let $A^*_{\mathbf{T}}$ denote the $\T$-equivariant Chow ring of point.
To the characters $t_i$,  we associate 
line bundles $L_i$ on $B\T$ and
Chern classes
$$s_i = c_1(L_i)\in A^*_{\mathbf{T}}$$
which generate the $\T$-equivariant Chow ring,
$$A^*_{\mathbf{T}} = \Z[s_1,s_2,s_3].$$
Let $(s_1,s_2,s_3)$ be the maximal ideal. Let
$$(A^*_{\mathbf{T}})_{loc} = \Q[s_1,s_2,s_3]_{(s_1,s_2,s_3)}$$
denote the localization.

Following the notation of Section \ref{ttt}, let 
$S^M_\mmuu$ be the set of components
of the
moduli space  of $\T$-invariant
submodules of $M/\O_{C_\mmuu}$.
Since the character \eqref{vertexchar} depends only
upon the local data at the vertex,
$\bV_{\mathcal{Q}}$ is well-defined for
$[\mathcal{Q}] \in S^M_\mmuu$. Let
\begin{eqnarray*}
\bw(\mathcal{Q}) & =&
\int_{\mathcal{Q}} e(T_\mathcal{Q}) \cdot e(-\bV_{\mathcal{Q}}) 
 \ \in (A^*_\T)_{loc}
\end{eqnarray*}
be the  integral of the evaluation of the
contribution \eqref{vertexchar} on $\mathcal{Q}$.
The integral is well-defined by Section \ref{fxy}.

Let $\ell(\mathcal{Q})$ denote the number of boxes{\footnote{As usual,
unlabelled type $\III$ boxes count twice.}} in
the labelled configuration associated to $\mathcal{Q}$.
Let $|\mmuu|$ denote the renormalized volume{\footnote{
The renormalized volume $|\pi|$ is defined by
$$
|\pi| = \#
\left\{\pi \cap [0,\dots,N]^3 \right\}-
(N+1) \sum_1^3 |\mu^i| \,, \quad N\gg 0 \,.
$$
The renormalized volume is
independent of the
cut-off $N$ as long as $N$ is
sufficiently large. The number $|\pi|$
so defined may be negative.
}} 
of the
partition $\pi$ corresponding to $\I_{C_\mmuu}$.
Finally, the {\em stable pairs equivariant vertex} is defined by
\begin{equation}\label{vvpe}
\bW_\mmuu^P =
\sum_{[\mathcal{Q}]\in S^M_\mmuu}
 \bw(\mathcal{Q})\ q^{\ell(\mathcal{Q})+|\mmuu|}\ \in
\Q(s_1,s_2,s_3)(\!(q)\!).
\end{equation}
\subsection{Vertex correspondence}
The non-normalized DT equivariant vertex $\bW_\mmuu^{DT,nn}$
 is defined in parallel terms in Section 4 of \cite{MNOP2}.
Define the {\em normalized DT equivariant vertex} by
$$\bW_\mmuu^{DT} = \frac{{\bW}_\mmuu^{DT,nn}}
{{\bW}_{\emptyset,\emptyset,\emptyset}^{DT,nn}}\ .$$
The degree 0 series has been calculated in \cite{MNOP2},
$$\bW_{\emptyset, \emptyset, \emptyset}^{DT,nn}= 
M(-q)^{- \frac{(s_1+s_2)(s_1+s_3)(s_2+s_3)}{s_1s_2s_3}}\,,$$
where $M$ denotes the MacMahon function
$$M(-q) = \prod_{n\geq 1} \frac{1}{(1-(-q)^n)^n}\,.$$

\begin{conj}\label{cqgt} The equivariant vertices agree,
$$\bW^P_{\mmuu} = \bW^{DT}_\mmuu.$$
\end{conj}

In fact, there is a straightforward approach to proving
Conjecture \ref{cqgt} following the path{\footnote{The
DT path in turn follows the Gromov-Witten trail \cite{BryanP,gwan}.}} 
already taken in
DT theory  \cite{dtan,moop,lcdt}. Each DT step can very likely
be followed by the identical step in the theory of stable pairs.
If the path is followed to the end, a proof of Conjecture
\ref{cqgt} will be obtained. The road is long and, even in
DT theory, missing foundation developments for the
relative geometry. One could hope for a more direct
proof of Conjecture \ref{cqgt} via a wall-crossing
analysis in the derived category.

\subsection{Example}
A basic example to consider is $\bW_{(1),\emptyset,\emptyset}^P$.
The calculation from the definition of the stable pairs
vertex is almost trivial.

\begin{lem} \label{xdd}
  $\bW_{(1),\emptyset,\emptyset}^P= (1+q)^{\frac{s_2+s_3}{s_1}}.$
\end{lem}
\begin{proof}
The outgoing partitions are $\mu^1=(1)$ and $\mu^2=\mu^3=\emptyset$.
The component set $S^M_{(1),\emptyset,\emptyset}$
is in bijective correspondence with the positive integers $k$.
The $\T$-fixed point $\mathcal{Q}^k$ corresponds to the length
$k$
box configuration in Figure 6.

\begin{figure} 
\center{\epsfig{file=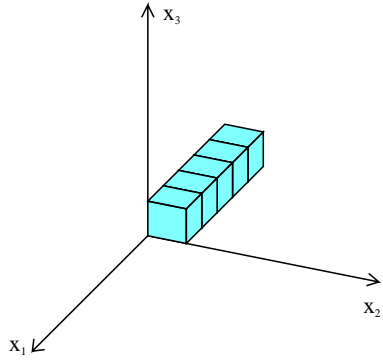, width=95mm}}
\vspace{-43mm}
\caption{The box configuration corresponding to $\mathcal{Q}^5$}
\end{figure}

The $q^k$ coefficient of $\bW_{(1),\emptyset,\emptyset}$ is
obtained by simply expanding the definition of 
$\bV_{\mathcal{Q}^k}$.
The $\T$-character $\FFF_{\mathcal{Q}^k}$ is
$$\FFF_{\mathcal{Q}^k} = \frac{t_1^{-k}}{1-t_1}$$
and the $\T$-character associated to the single edge $e$
along the $x_1$-axis is
$$\FFF_{e} = 1.$$
Unwinding \eqref{gx34} yields
$$\bV_{\mathcal{Q}^k} = \sum_{i=1}^k t_1^{-i} - \sum_{i=0}^{k-1} 
\frac{t_1^{i}}{t_2t_3}\,.$$
In fact, only the first two terms of \eqref{vertexchar}
contribute --- the others cancel with the redistribution.
Then,
\begin{eqnarray*}
\bw(\mathcal{Q}^k)& = &
\int_{\mathcal{Q}_k} e(-\bV_{\mathcal{Q}^k})\\
& = & \frac{(-s_2-s_3)(s_1-s_2-s_3) \cdots ((k-1)s_1-s_2-s_3)}
{(-s_1) (-2s_1) \cdots (-ks_1)} \\
& = & \frac{1}{k!} \left( \frac{s_2+s_3}{s_1}\right)
\left( \frac{s_2+s_3}{s_1}-1\right)\cdots
\left( \frac{s_2+s_3}{s_1}-(k-1)\right),
\end{eqnarray*}
which is clearly the $q^k$ coefficient of $(1+q)^{\frac{s_2+s_3}{s_1}}$.
\end{proof}

The result agrees with the DT  calculation
of $\bW^{DT}_{(1),\emptyset,\emptyset}$
 in Section 6 of \cite{lcdt} and verifies Conjecture
\ref{cqgt}. The intricacy of the DT argument using 
localization relations, divisibility, and the Hilbert-Chow
morphism shows the difference in the theories.

\section{Toric Calabi-Yau vertex} \label{cyver}
\subsection{Calabi-Yau torus}
Let $X$ be a toric Calabi-Yau 3-fold with canonical form $\omega$.
The main example is a {\em local Calabi-Yau toric surface},
the total space of the canonical bundle
$$K_S \rightarrow S$$
of toric surface $S$.
As in Section \ref{tocy}, let $\T_0\subset T$ be the
torus preserving $\omega$. The stable pairs invariants
take a much simpler form when computed $\T_0$-equivariantly.

The direct approach to calculating the $\T_0$-equivariant
vertex is simply by restriction,
$$\bW^P_{\mmuu,CY} = \bW^P_{\mmuu} \ |\_{s_1+s_2+s_3=0}\,.$$
The main drawback is the summation
over $[\mathcal{Q}]\in S^M_\mmuu$ in the
definition \eqref{vvpe}
of $\bW^P_\mmuu$ has to be done
{\em before} restriction as the individual summands
may have poles at $s_1+s_2+s_3=0$.
A better approach is to analyse the $\T_0$-localization
formula geometrically.

\subsection{$\T_0$-localization}
Let $S^M_{\mmuu,0}$ be the set of connected components
 of the moduli space
of $\T_0$-invariant submodules of $M/O_{C_\mmuu}$.
By Conjecture 2, such components
$\mathcal{Q}_0$ are nonsingular.{\footnote{We
assume the validity of Conjecture \ref{vvt} throughout the Section.}}
The quotient torus $$\C^* \cong \T/\T_0$$
 acts on each $\mathcal{Q}_0$.
The $\C^*$-fixed loci of $\mathcal{Q}_0$ are
elements of $[\mathcal{Q}]\in S^M_\mmuu$.
Conversely, every element of $[\mathcal{Q}]\in S^M_\mmuu$ arises
as a $\C^*$-fixed locus of a unique
$[\mathcal{Q}_0]\in S^M_{\mmuu,0}$.

The scheme $\mathcal{Q}_0$ carries a $\T_0$-fixed
obstruction theory,
\begin{equation}\label{bbyy}
(E^{-1}_{\mathcal{Q}_0})^{\T_0}
 \rightarrow (E^{0}_{\mathcal{Q}_0})^{\T_0}.
\end{equation}
By the nonsingularity of $\mathcal{Q}_0$, the virtual
class $[\mathcal{Q}_0]$ is the the Euler class of
an obstruction bundle. By the self-duality of the
obstruction theory \eqref{bbyy},
$$[\mathcal{Q}_0]^{vir} = e(\Omega_{\mathcal{Q}_0}) \cap [\mathcal{Q}_0].$$
Then,
by the virtual localization formula for the $\C^*$-action and self-duality again,
\begin{eqnarray*}
\bW^P_{\mmuu} \ |_{s_1+s_2+s_3=0} &  = &
\sum_{[\mathcal{Q}_0]\in S^M_{\mmuu,0}} \int_{\mathcal{Q}_0}
e(\Omega_{\mathcal{Q}_0})
\ 
(-1)^{\text{rk}(E^0_{\mathcal{Q}_0})^{m_0}}\\
&= &
\sum_{[\mathcal{Q}_0]\in S^M_{\mmuu,0}} \chi\_{top}(\mathcal{Q}_0)
\ (-1)^{\text{dim}_\C \mathcal{Q}_0}
(-1)^{\text{rk}(E^0_{\mathcal{Q}_0})^{m_0}},
\end{eqnarray*}
where
$\text{rk}(E^0_{\mathcal{Q}_0})^{m_0}$ is the rank of the
summand with nontrivial $\T_0$-weight.

\subsection{Theorem \ref{ddd}}
Using the Euler characteristic identity
$$\chi\_{top}(\mathcal{Q}_0) = \sum_{[\mathcal{Q}]\in S^M_{\mmuu}, \ 
\mathcal{Q} \subset \mathcal{Q}_0} \chi\_{top}(\mathcal{Q}),$$
we can rewrite the Calabi-Yau vertex as
\begin{multline}\label{xbp}
\bW^P_{\mmuu} \ |_{s_1+s_2+s_3=0} =\\
\sum_{[\mathcal{Q}_0]\in S^M_{\mmuu,0}} \
\sum_{[\mathcal{Q}]\in S^M_{\mmuu}, \ \mathcal{Q} \subset \mathcal{Q}_0} 
\chi\_{top}(\mathcal{Q})
\ (-1)^{\text{dim}_\C \mathcal{Q}_0}
(-1)^{\text{rk}(E^0_{\mathcal{Q}})^{m_0}}.
\end{multline}

Fix an element $[Q]\in\mathcal{Q} \subset \mathcal{Q}_0$.
The calculation of $$(-1)^{\text{dim}_\C \mathcal{Q}_0}
(-1)^{\text{rk}(E^0_{Q})^{m_0}}$$
precisely follows Sections 4.10-4.11 of \cite{MNOP1}. 
Let
$$\bV_Q(t_1,t_2,t_3)= \bV_{\mathcal{Q}} |_{[Q]},$$
where the restriction to the point $[Q]\in \mathcal{Q}$ 
kills all the $K$-theory information.
The strategy of \cite{MNOP1} is to split the vertex contribution
as
$$\bV_Q(t_1,t_2,t_3)=
\bV_Q^+ + \bV_Q^-$$
where
$${\bV}_Q^+ = -\bV_Q^-.$$
We use precisely the same formulas to define the splitting
here, see Section 4.11 of \cite{MNOP1}.
Then,
$$\bV^+_Q(1,1,1)= \mathsf{Con}(\bV^+_Q|_{t_1t_2t_3=1}) +
\text{rank}(E^0_{Q})^{m_0} \ \mod 2.$$
The constant term is treated exactly as
in Section 4.11 of \cite{MNOP1}. The conclusion, using
Conjecture \ref{vvgg}, is
$$\mathsf{Con}(\bV^+_Q|_{t_1t_2t_3=1}) = {\text{dim}_\C \mathcal{Q}_0}.$$
Finally, from the definition of the splitting,
$$\bV^+_Q(1,1,1)= \ell(Q)+ |\mmuu|,$$
where the second summand is the renormalized volume.
We conclude
\begin{equation}\label{blww}
\bW^P_{\mmuu} \ |_{s_1+s_2+s_3=0} =
\sum_{[\mathcal{Q}]\in S^M_{\mmuu}} 
\chi\_{top}(\mathcal{Q}) \ (-q)^{\ell(\mathcal{Q})+ |\mmuu|}\,.
\end{equation}

Our calculation of the $\T_0$-vertex \eqref{blww}
is complete in the 1 and 2-leg cases and conjectural in
the 3-leg case (since Conjectures \ref{vvt} and \ref{vvgg} are
unproven there). In particular, the result is established
for local Calabi-Yau toric surfaces.
Though we have not found examples violating
Conjecture \ref{vvt}, perhaps \eqref{blww} holds
without requiring the nonsingularity of the components $\mathcal{Q}_0$.

The statement of Theorem \ref{ddd} follows formally
from the vertex equation \eqref{blww} and the edge calculus
of Section 4.10 of \cite{MNOP1}. The edge calculus is
identical for stable pairs and ideal sheaves.\qed

\subsection{Example}
The example of the Calabi-Yau vertex  with outgoing partitions
$$\mu^1=(1),\  \mu^2=(1), \ \mu^3=(1)$$
was worked out at the end of \cite{pt}. 

We consider here the vertex $\bW^P_{(1),(2),(1),CY}$
with outgoing partitions
$$\mu^1=(1),\ \mu^2=(2),\  \mu^3=(1).$$
The scheme-theoretic support curve $C_{(1),(2),(1)}$
is pictured in Figure 7
with renormalized volume  $-3$.
\begin{figure} 
\center{\epsfig{file=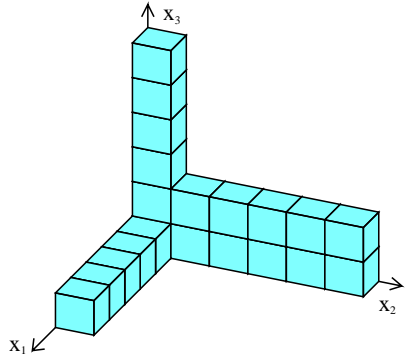, width=95mm}}
\vspace{-43mm}
\caption{The $\T$-character of the support curve $C_{(1),(2),(1)}$.}
\end{figure}
The module $M /\O_{C_{\mmuu}}$ is pictured in
Figure 8. There is a single type $\II$ box at $x_3$ and
a single type $\III$ box at $1$, over the origin.
\begin{figure} 
\center{\epsfig{file=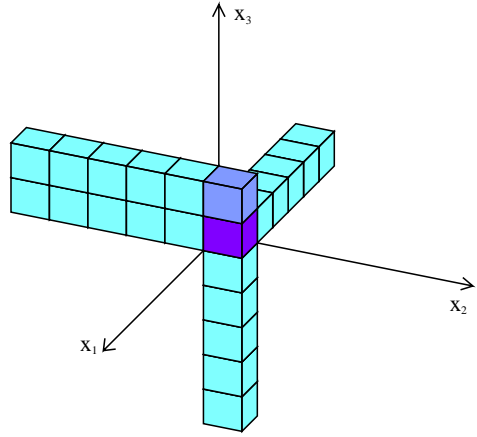, width=95mm}}
\vspace{-43mm}
\caption{The $\T$-character of the module $M /\O_{C_{(1),(2),(1)}}$.}
\end{figure}

We now count box configurations up to length 3. There is a unique
empty box configuration corresponding to the 0 submodule.
There are 2 distinct configurations of length 1:
\begin{enumerate}
\item[(i)]
a box at $1$ labelled by $\C\cdot \mathbf{1}_{0}\in\PP^1$,
\item[(ii)] a box supported at $x_3$.
\end{enumerate}
The box  (i) 
is annihilated by multiplication by all the  $x_i$ in the module
structure.
There are 3 types of configurations of length 2:
\begin{enumerate}
\item[(i)]  boxes at
$x_2^{-1} x_3$ and $x_3$,
\item[(ii)] a box at 1 labelled with $\C\cdot \mathbf{1}_0$ and a box
at $x_1^{-1}$, 
\item[(iii)] a box at 1 with any label in $\PP^1$ and a box at $x_3$.
\end{enumerate}
The moduli space in (iii) is $\PP^1$.
The configurations of length 3 are classified as follows:
\begin{enumerate}
\item[(i)] an unlabelled (length 2) box at $1$ and a box at $x_3$,
\item[(ii)]
 a box at 1 with any label in $\PP^1$  and boxes at $x_2^{-1}x_3$ and $x_3$,
\item[(iii)] a box at $x_1^{-1}$, a box at 1 labelled by $\C\cdot \mathbf{1}_0$,
and a box at $x_3$, 
\item[(iv)] a box at
$x_3^{-1}$,  a box at 1 labelled by $\C\cdot \mathbf{3}_0$, and
a box at $x_3$,
\item[(v)] a box at 1 labelled by $\C\cdot \mathbf{1}_0$ and
boxes at $x_1^{-2}$ and $x_1^{-1}$,
\item[(vi)] boxes at $x_2^{-2}x_3$, $x_2^{-1}x_3$, and $x_3$.
\end{enumerate}
The moduli space in (ii) is $\PP^1$.

Equation \eqref{blww} for the Calabi-Yau vertex and
the above box counting (including the
Euler characteristics) together yield {\footnote{
A closed formula for the Calabi-Yau
vertex here can easily be found.}}
$$(-q)^3\ \bW^P_{(1),(2),(1),CY} = 1+2(-q)+4(-q)^2+7(-q)^3+\ldots\ .
$$
Standard box counting yields the non-normalized
DT vertex
$$
(-q)^3\ \bW^{DT,nn}_{(1),(2),(1),CY}
=1+3(-q)+9(-q)^2+23(-q)^3+\ldots\ ,
$$
see \cite{MNOP1}.
Using the MacMahon series
$$M(-q)= 1+(-q)+3 (-q)^2+6(-q)^3 +\ldots\ ,$$
we can check Conjecture \ref{cqgt}
up to order 3,
$$1+2(-q)+4(-q)^2+7(-q)^3+\ldots = \frac{
1+3(-q)+9(-q)^2+23(-q)^3+\ldots}
{1+(-q)+3(-q)^2+6(-q)^3+\ldots}\,.
$$

\section{Descendents}\label{desc}
\subsection{Chern characters}
Let $X$ be a nonsingular
toric 3-fold. Consider the
$\T$-equivariant descendent invariants
\begin{multline}\label{vpzz}
\left\langle \tau_{i_1}(\gamma_1) \cdots \tau_{i_k}(\gamma_k) \right
\rangle_{n,\beta}^X =\\
\int_{P_n(X,\beta)} \prod_{k=1}^m \tau_{i_k}(\gamma_k)
\Big( [P_{n}(X,\beta)]^{vir}\Big)\in A^*_\T,
\end{multline}
where $\gamma_i \in A^*_\T(X,\mathbb{Z})$.
The operators $\tau_i(\gamma)$ are defined by
$$
\pi_{P*}\big(\pi_X^*(\gamma)\cdot \text{ch}_{2+i}(\FF)
\cap(\pi_P^*(\ \cdot\ )\big)\colon 
A^\T_*(P_{n}(X,\beta))\to A^\T_*(P_{n}(X,\beta))
$$
where we follow the notation of  Section \ref{des}.

In order to calculate \eqref{vpzz} by $\T$-localization, we
must determine the action of
the operators $\tau_{i}(\gamma)$ on the $\T$-equivariant
cohomology of the $\T$-fixed loci 
$\mathbf{Q}$.

\subsection{Local to global}
Let $\mathbf{Q}$ be a component of the $\T$-fixed locus of 
$P_n(X,\beta)$ determined by the local data 
$$\{\mathcal{Q}_\alpha\}_{[X_\alpha]\in \mathsf{V}(X)}, \ \ 
\{\mu_{\alpha\beta}\}_{[C_{\alpha\beta}]\in\mathsf{E}(X)}\,.$$
On the chart $U_\alpha\subset X$, the Chern character
$$\text{ch}_l(\FF_\alpha)|_{\mathbf{Q}_\alpha\times [X_\alpha]} \in
K(\mathcal{Q}_\alpha) \otimes_\Z \Z[t_1^\pm,t_2^\pm,t_3^\pm]
$$ 
 is determined completely by the character $\FFF_\alpha$,
\begin{eqnarray}\label{ftrr}
\text{ch}_l(\FF_\alpha)|_{\mathbf{Q}_\alpha\times [X_\alpha]}  &=&
\text{ch}_l(1+P_\alpha) \\ \nonumber
& = & \text{ch}_l\big(\FFF_\alpha\cdot (1-t_1)(1-t_2)(1-t_3)\big),
\end{eqnarray}
where we follow the notation of Section \ref{vc}. In particular,
the second equality is \eqref{PQ}.
For notational convenience, we denote 
\eqref{ftrr} and the pull-back to $\mathbf{Q}$ simply by
$\text{ch}_l(\FF_\alpha)$.

The contribution of $\mathbf{Q}$ to the descendent
\eqref{vpzz} is calculated by summing over
all distributions
$$\sigma: j \mapsto [X_{\sigma(j)}]\in \mathsf{V}(X)$$ of
 the
insertions $\tau_{i_j}(\gamma_j)$ to the vertices
and integrating
\begin{multline*}
\sum_\sigma \int_{\mathbf{Q}}
e(T_{\mathbf{Q}}) \cdot e\left(-\sum 
\mathsf{V}_\alpha
-\sum \mathsf{E}_{\alpha\beta}\right)
\cdot \prod_{j=1}^k \text{ch}_{2+i_j}(\FF_{\sigma(j)}) 
\frac{\gamma_j|_{[X_{\sigma(j)}]}}{e(T_{\sigma(j)})}
\end{multline*}
where $T_{\sigma(j)}$ is the $\T$-equivariant tangent space
to $X_{\sigma(j)}$.

\subsection{Descendent vertex}
The {\em descendent vertex} $\bW_\mmuu^p(\tau_{i_1} \cdots \tau_{i_k})$
is obtained from the descendent weight
\begin{multline}
\bw_{\tau_{i_1}\cdots \tau_{i_k}}
(\mathcal{Q}) =\\
\int_{\mathcal{Q}} e(T_\mathcal{Q}) \cdot e(-\bV_{\mathcal{Q}})
\cdot \prod_{j=1}^k  
\text{ch}_{2+i_j}\big(\FFF_{\mathcal{Q}}\cdot (1-t_1)(1-t_2)(1-t_3)\big)
\end{multline}
taking values in $(A^*_\T)_{loc}$.
By definition,
\begin{equation}\label{vvped}
\bW_\mmuu^P(\tau_{i_1} \cdots \tau_{i_k})
 =
\sum_{[\mathcal{Q}]\in S^M_\mmuu}
 \bw_{\tau_{i_1} \cdots \tau_{i_k}}
(\mathcal{Q})\ q^{\ell(\mathcal{Q})+|\mmuu|}\ \in
\Q(s_1,s_2,s_3)(\!(q)\!).
\end{equation}
\label{heyle}

\subsection{Example}
As an example, we calculate the descendent vertex 
$$\bW^P_{(1),\emptyset,\emptyset}(\tau_i)$$ from the definitions.

Following the notation of Lemma \ref{xdd}, we have
$$
\FFF_{\mathcal{Q}^k}\cdot (1-t_1)(1-t_2)(1-t_3) =  
t_1^{-k}(1-t_2)(1-t_3).
$$
Then, by Lemma \ref{xdd} and \eqref{vvped},
\begin{eqnarray*}
\sum_{i\geq -2}\bW^P_{(1),\emptyset,\emptyset}(\tau_i)z^{2+i}& = &
\sum_{k\geq 0} \sum_{i\geq -2} \bw_{\tau_{i}}(\mathcal{Q}^k)z^{2+i}\\
& =& \nonumber
(1-e^{zs_2})(1-e^{zs_3})
\exp(-zs_1 q\frac{d}{dq})\left[ (1+q)^{\frac{s_2+s_3}{s_1}}\right].
\end{eqnarray*}
The $\tau_{-2}$ and $\tau_{-1}$ terms are included formally.
The formulas for
$$\bW^P_{(1),\emptyset,\emptyset}(\tau_{i_1}\cdots \tau_{i_k})$$
are no more difficult.

Let $[L]\in H_2(\PP^3,\Z)$ be the class of a line.
We can use the determination of $\bW^P_{(1),\emptyset,\emptyset}(\tau_5)$
to calculate the $\T$-equivariant descendent series of $\PP^3$
in degree 1,
\begin{equation}\label{ddssq}
Z^{\PP^3}_{P,[L]}\big(   \prod_{j=1}^k \tau_{i}(1)
\big) \ \in \Q[s_1,s_2,s_3](\!(q)\!).
\end{equation}
Expansion of the $\T$-contribution formula of Section \ref{heyle}
then immediately yields a qualitative result: the descendent
series \eqref{ddssq}
are {\em rational} in $q$,
\begin{equation*}
Z^{\PP^3}_{P,[L]}\big(   \prod_{j=1}^k \tau_{i}(1)
\big) \ \in \Q[s_1,s_2,s_3]\otimes_\Q \Q(q).
\end{equation*}

In the non-equivariant limit, the rationality provides evidence
for Conjecture \ref{111}. However, we speculate, in the toric case, that 
rationality holds for all such  $\T$-equivariant
descendent series.

 The descendent invariants for the
theory of stable pairs appear much better behaved than the
descendents in DT theory. Because of the wandering points, the
DT descendent series corresponding to \eqref{ddssq} --- even in
the non-equivariant case --- contain
significant irrationalities. The Gromov-Witten descendent
series also contain irrationalities. The framework of
the GW/DT correspondence for descendents was outlined in \cite{MNOP2}.
The correspondence with the 
descendent theory of stable pairs remains to be fully explored.

\vspace{+8 pt}
\noindent
Department of Mathematics\\
Princeton University\\
rahulp@math.princeton.edu

\vspace{+8 pt}
\noindent
Department of Mathematics \\
Imperial College \\
rpwt@imperial.ac.uk

\end{document}

%% file: boxes.pstex_t
\begin{picture}(0,0)%
\includegraphics{boxes.pstex}%
\end{picture}%
\setlength{\unitlength}{1973sp}%
\begingroup\makeatletter\ifx\SetFigFont\undefined%
\gdef\SetFigFont#1#2#3#4#5{%
  \reset@font\fontsize{#1}{#2pt}%
  \fontfamily{#3}\fontseries{#4}\fontshape{#5}%
  \selectfont}%
\fi\endgroup%
\begin{picture}(3848,3147)(1801,-4990)
\put(2176,-2760){\makebox(0,0)[lb]{\smash{{\SetFigFont{6}{7.2}{\rmdefault}{\mddefault}{\updefault}{\color[rgb]{0,0,0}$x_2^3$}%
}}}}
\put(2720,-2760){\makebox(0,0)[lb]{\smash{{\SetFigFont{6}{7.2}{\rmdefault}{\mddefault}{\updefault}{\color[rgb]{0,0,0}$x_1x_2^3$}%
}}}}
\put(2720,-3361){\makebox(0,0)[lb]{\smash{{\SetFigFont{6}{7.2}{\rmdefault}{\mddefault}{\updefault}{\color[rgb]{0,0,0}$x_1x_2^2$}%
}}}}
\put(2176,-3361){\makebox(0,0)[lb]{\smash{{\SetFigFont{6}{7.2}{\rmdefault}{\mddefault}{\updefault}{\color[rgb]{0,0,0}$x_2^2$}%
}}}}
\put(2720,-3961){\makebox(0,0)[lb]{\smash{{\SetFigFont{6}{7.2}{\rmdefault}{\mddefault}{\updefault}{\color[rgb]{0,0,0}$x_1x_2$}%
}}}}
\put(2176,-3961){\makebox(0,0)[lb]{\smash{{\SetFigFont{6}{7.2}{\rmdefault}{\mddefault}{\updefault}{\color[rgb]{0,0,0}$x_2$}%
}}}}
\put(2176,-4561){\makebox(0,0)[lb]{\smash{{\SetFigFont{6}{7.2}{\rmdefault}{\mddefault}{\updefault}{\color[rgb]{0,0,0}$1$}%
}}}}
\put(2720,-4561){\makebox(0,0)[lb]{\smash{{\SetFigFont{6}{7.2}{\rmdefault}{\mddefault}{\updefault}{\color[rgb]{0,0,0}$x_1$}%
}}}}
\put(3376,-4561){\makebox(0,0)[lb]{\smash{{\SetFigFont{6}{7.2}{\rmdefault}{\mddefault}{\updefault}{\color[rgb]{0,0,0}$x_1^2$}%
}}}}
\put(3976,-4561){\makebox(0,0)[lb]{\smash{{\SetFigFont{6}{7.2}{\rmdefault}{\mddefault}{\updefault}{\color[rgb]{0,0,0}$x_1^3$}%
}}}}
\put(4576,-4561){\makebox(0,0)[lb]{\smash{{\SetFigFont{6}{7.2}{\rmdefault}{\mddefault}{\updefault}{\color[rgb]{0,0,0}$x_1^4$}%
}}}}
\put(3337,-3961){\makebox(0,0)[lb]{\smash{{\SetFigFont{6}{7.2}{\rmdefault}{\mddefault}{\updefault}{\color[rgb]{0,0,0}$x_1^2x_2$}%
}}}}
\put(3313,-3361){\makebox(0,0)[lb]{\smash{{\SetFigFont{6}{7.2}{\rmdefault}{\mddefault}{\updefault}{\color[rgb]{0,0,0}$x_1^2x_2^2$}%
}}}}
\put(5326,-4900){\makebox(0,0)[lb]{\smash{{\SetFigFont{7}{8.4}{\rmdefault}{\mddefault}{\updefault}{\color[rgb]{0,0,0}$x_1$}%
}}}}
\put(1750,-2011){\makebox(0,0)[lb]{\smash{{\SetFigFont{7}{8.4}{\rmdefault}{\mddefault}{\updefault}{\color[rgb]{0,0,0}$x_2$}%
}}}}
\end{picture}%

%% file: moves.pstex_t
\begin{picture}(0,0)%
\includegraphics{moves.pstex}%
\end{picture}%
\setlength{\unitlength}{3158sp}%
\begingroup\makeatletter\ifx\SetFigFont\undefined%
\gdef\SetFigFont#1#2#3#4#5{%
  \reset@font\fontsize{#1}{#2pt}%
  \fontfamily{#3}\fontseries{#4}\fontshape{#5}%
  \selectfont}%
\fi\endgroup%
\begin{picture}(5236,3505)(1801,-5519)
\put(5476,-5461){\makebox(0,0)[lb]{\smash{{\SetFigFont{10}{12.0}{\rmdefault}{\mddefault}{\updefault}{\color[rgb]{0,0,0}$x_2$}%
}}}}
\put(1801,-2761){\makebox(0,0)[lb]{\smash{{\SetFigFont{10}{12.0}{\rmdefault}{\mddefault}{\updefault}{\color[rgb]{0,0,0}$x_3$}%
}}}}
\put(4525,-5161){\makebox(0,0)[lb]{\smash{{\SetFigFont{5}{6.0}{\rmdefault}{\mddefault}{\updefault}{\color[rgb]{0,0,0}$x_2^{a_2\!+\!A}$}%
}}}}
\put(2713,-3961){\makebox(0,0)[lb]{\smash{{\SetFigFont{5}{6.0}{\rmdefault}{\mddefault}{\updefault}{\color[rgb]{0,0,0}$x_2^{a_2}x_3^{a_3}$}%
}}}}
\put(3913,-4561){\makebox(0,0)[lb]{\smash{{\SetFigFont{5}{6.0}{\rmdefault}{\mddefault}{\updefault}{\color[rgb]{0,0,0}$x_2^{\!a_2\!+\!m}\!\!x_3^{a_3\!-\!n}$}%
}}}}
\put(3913,-3961){\makebox(0,0)[lb]{\smash{{\SetFigFont{5}{6.0}{\rmdefault}{\mddefault}{\updefault}{\color[rgb]{0,0,0}$x_2^{\!a_2\!+\!m}\!x_3^{a_3}$}%
}}}}
\put(6376,-3361){\makebox(0,0)[lb]{\smash{{\SetFigFont{5}{6.0}{\rmdefault}{\mddefault}{\updefault}{\color[rgb]{0,0,0}$x_2^{a_2+A-w}\!x_3^{-w}$}%
}}}}
\put(3976,-2086){\makebox(0,0)[lb]{\smash{{\SetFigFont{5}{6.0}{\rmdefault}{\mddefault}{\updefault}{\color[rgb]{0,0,0}$x_2^{a_2-w}\!x_3^{a_3-w}$}%
}}}}
\put(5326,-2086){\makebox(0,0)[lb]{\smash{{\SetFigFont{5}{6.0}{\rmdefault}{\mddefault}{\updefault}{\color[rgb]{0,0,0}$x_2^{a_2+m-w}\!x_3^{a_3-w}$}%
}}}}
\put(5851,-2686){\makebox(0,0)[lb]{\smash{{\SetFigFont{5}{6.0}{\rmdefault}{\mddefault}{\updefault}{\color[rgb]{0,0,0}$x_2^{a_2+m-w}\!x_3^{a_3-n-w}$}%
}}}}
\put(5101,-3136){\makebox(0,0)[lb]{\smash{{\SetFigFont{10}{12.0}{\rmdefault}{\mddefault}{\updefault}{\color[rgb]{0,0,0}$\phi$}%
}}}}
\end{picture}%